\documentclass[11pt,a4paper,twoside,final]{scrartcl}
\usepackage{a4wide}
\usepackage{amsfonts}
\usepackage{amsmath}
\usepackage{amssymb}
\usepackage{amsthm}
\usepackage{mathtools}
\usepackage[mathscr]{euscript}
\usepackage[numbers]{natbib}
\usepackage{dsfont}
\usepackage{booktabs}
\usepackage[Algorithm]{algorithm}
\usepackage{algorithmic}
\usepackage{pgfplots}
\usepackage{pgfplotstable}
\usepackage{graphicx}
\usepackage{soul}
\usepackage{color}
\usepackage{colortbl}
\usepackage{multirow}
\usepackage{boxedminipage}
\usepackage{enumerate}
\usepackage{collcell}
\usepackage{caption}
\usepackage{subcaption}
\usepackage{rotating}
\usepackage{hyperref}

\newcommand{\A}{\ensuremath{\mathcal{A}}}
\newcommand{\C}{\ensuremath{\mathbb{C}}}
\newcommand{\D}{\ensuremath{\mathcal{D}}}

\newcommand{\F}{\ensuremath{\mathcal{F}}}
\newcommand{\G}{\ensuremath{\mathcal{G}}}
\renewcommand{\H}{\ensuremath{\mathcal{H}}}
\newcommand{\I}{\ensuremath{\mathrm{I}}}

\renewcommand{\L}{\ensuremath{\mathcal{L}}}
\newcommand{\N}{\ensuremath{\mathbb{N}}}

\newcommand{\R}{\ensuremath{\mathbb{R}}}

\newcommand{\T}{\ensuremath{\mathbb{T}}}

\newcommand{\U}{\ensuremath{\mathcal{U}}}

\newcommand{\bolda}{{\ensuremath{\boldsymbol{a}}}}

\newcommand{\boldb}{{\ensuremath{\boldsymbol{b}}}}

\newcommand{\boldf}{{\ensuremath{\boldsymbol{f}}}}

\newcommand{\boldk}{{\ensuremath{\boldsymbol{k}}}}

\newcommand{\boldell}{{\ensuremath{\boldsymbol{\ell}}}}

\newcommand{\boldx}{{\ensuremath{\boldsymbol{x}}}}
\newcommand{\boldy}{{\ensuremath{\boldsymbol{y}}}}
\newcommand{\boldz}{{\ensuremath{\boldsymbol{z}}}}

\newcommand{\boldtheta}{{\ensuremath{\boldsymbol{\theta}}}}
\newcommand{\boldxi}{{\ensuremath{\boldsymbol{\xi}}}}

\newcommand{\boldone}{{\ensuremath{\boldsymbol{1}}}}

\newcommand{\e}{\textnormal{e}}
\newcommand{\ii}{\textnormal{i}}

\newcommand{\ghat}{{\ensuremath{\hat{g}}}}

\newcommand{\abs}[1]{\left\vert #1\right\vert}
\newcommand{\norm}[1]{\left\Vert #1\right\Vert}

\newcommand{\diff}{\ensuremath{\mathrm{d}}}

\definecolor{color1}{HTML}{E69F00}
\definecolor{color2}{HTML}{56B4E9}
\definecolor{color3}{HTML}{009E73}
\definecolor{color4}{HTML}{F0E442}
\definecolor{color5}{HTML}{0072B2}
\definecolor{color6}{HTML}{D55E00}
\definecolor{color7}{HTML}{CC79A7}

\newtheorem{theorem}{Theorem}[section]
\newtheorem{lemma}[theorem]{Lemma}
\newtheorem{remark}[theorem]{Remark}
\newtheorem{generalisation}[theorem]{Generalisation}
\newtheorem{definition}[theorem]{Definition}
\newtheorem{example}[theorem]{Example}
\newtheorem{corollary}[theorem]{Corollary}
\newtheorem{proposition}[theorem]{Proposition}

\usepgfplotslibrary{statistics}
\usetikzlibrary{patterns}
\pgfplotsset{compat=newest}

\pgfplotsset{
        my boxplot style/.style={
            boxplot prepared,
            draw=black,
            solid,
            fill=white,
            mark=*,
            every mark/.append style={
                fill=gray,
            },
        },
    }

\pgfplotsset{
    boxplot prepared from table/.code={
        \def\tikz@plot@handler{\pgfplotsplothandlerboxplotprepared}%
        \pgfplotsset{
            /pgfplots/boxplot prepared from table/.cd,
            #1,
        }
    },
    /pgfplots/boxplot prepared from table/.cd,
        table/.code={\pgfplotstablecopy{#1}\to\boxplot@datatable},
        row/.initial=0,
        make style readable from table/.style={
            #1/.code={
                \pgfplotstablegetelem{\pgfkeysvalueof{/pgfplots/boxplot prepared from table/row}}{##1}\of\boxplot@datatable
                \pgfplotsset{boxplot/#1/.expand once={\pgfplotsretval}}
            }
        },
        make style readable from table=lower whisker,
        make style readable from table=upper whisker,
        make style readable from table=lower quartile,
        make style readable from table=upper quartile,
        make style readable from table=median,
        make style readable from table=lower notch,
        make style readable from table=upper notch,
        make style readable from table=draw position,
        make style readable from table=box extend,
}
\makeatother

\numberwithin{equation}{section}
\numberwithin{table}{section}
\numberwithin{figure}{section}

\title{An approach to discrete operator learning based on sparse high-dimensional approximation}
\date{\today}
\author{Daniel Potts\footnotemark[1], Fabian Taubert\footnotemark[2]}

\hypersetup{pdfauthor=...,
          pdftitle={An approach to discrete operator learning based on sparse high-dimensional approximation},
          pdfsubject={},
          pdfcreator = {pdflatex},
          plainpages=false,
          pdfstartview=FitH,       
          pdfview=FitH,            
          pdfpagemode=UseOutlines, 
          bookmarksnumbered=true, 
          bookmarksopen=false,     
          bookmarksopenlevel=0,   
          colorlinks=true,       
          linkcolor=black,
          citecolor=black,
          urlcolor=black}

\begin{document}

\maketitle

\begin{abstract}
\small
We present a dimension-incremental method for function approximation in bounded orthonormal product bases to learn the solutions of various differential equations. Therefore, we decompose the source function of the differential equation into parameters like Fourier or Spline coefficients and treat the solution of the differential equation as a high-dimensional function w.r.t.\ the spatial variables, these parameters and also further possible parameters from the differential equation itself. Finally, we learn this function in the sense of sparse approximation in a suitable function space by detecting coefficients of the basis expansion with the largest absolute values. Investigating the corresponding indices of the basis coefficients yields further insights on the structure of the solution as well as its dependency on the parameters and their interactions and allows for a reasonable generalization to even higher dimensions and therefore better resolutions of the decomposed source function.

\medskip

\end{abstract}

\footnotetext[1]{
  Chemnitz University of Technology, Faculty of Mathematics, 09107 Chemnitz, Germany\\
  potts@mathematik.tu-chemnitz.de
}
\footnotetext[2]{
  Chemnitz University of Technology, Faculty of Mathematics, 09107 Chemnitz, Germany\\
  fabian.taubert@mathematik.tu-chemnitz.de
}
\medskip


\section{Introduction}
\label{sec:introduction}

In mathematical analysis, partial differential equations (PDEs) stand as formidable tools when it comes to modeling diverse phenomena across various scientific disciplines from fluid dynamics to quantum mechanics. Unfortunately, solving PDEs analytically as well as numerically can be quite difficult and thus became a challenging task. The numerical solution of PDEs is investigated thoroughly already since the mid-20th century until today by various well-known methods like finite difference methods \cite{Leveque07}, spectral methods \cite{Trefethen2000,CaHe23} or finite element methods (FEMs) \cite{ZiTa05}. On the other hand, with the age of artificial intelligence (AI), several new methods using machine learning techniques are currently arising and investigated for the solution of PDEs, including physics-informed neural networks (PINNs) \cite{RaPeKa19,Ka21,XiYuWa23,LuZh23}, convolutional neural networks \cite{GrHeKl24}, deep operator networks \cite{GoYiYuKa22, DeTeGi23, LiMoZh23} multilevel Picard approximations \cite{WeHuJeKr19,HuNg22} and neural operators \cite{Li21,JiMeLu22,Ni23,LiMiPeKaMi23,LaStTr24}, among numerous others.

It is not necessary to emphasize that both the classic and the AI methods have various advantages and disadvantages, which are also getting studied frequently in the last years \cite{BlEr21,GrKoLaSch23}. Especially for high-dimensional PDEs and for many-query settings, where multiple solutions of the PDE w.r.t.\ varying parameters, initial or boundary conditions are needed, machine learning algorithms have proven to outperform classical methods significantly. One of the main reasons for this is the great performance of neural networks in the framework of operator learning on PDEs, i.e., learning the underlying solution operator of a PDE which maps initial and/or boundary conditions as well as other parameters to the PDE solution. For a comprehensive overview of operator learning theory, algorithms, and applications, we refer the reader to the recent surveys \cite{BoTo24,KoLaStu24}.

However, for operator learning of simpler differential equations, classical methods can still compete with machine learning. As an example, consider the one-dimensional differential equation 
\begin{align}\label{eq:intro_pde}
\L u = f 
\end{align}
with the differential operator $\L = \frac{d}{dx}$ and some initial condition $u(0) = 0$. Assume the right-hand side function $f$ to be given (or at least be well approximated) by a partial Fourier sum
\begin{align*}
f(x) = \sum_{\ell=0}^{N-1} f_\ell \e^{2\pi\ii \ell x}.
\end{align*}
If we denote the vector of Fourier coefficients $\boldf = (f_0,\ldots,f_{N-1}) \in \C^N$, we now aim for a solution $u$ of the form 
\begin{align} \label{eq:intro_sol}
u(x,\boldf) = \sum_{\boldk \in I} u_\boldk T_\boldk(x,\boldf)
\end{align}
with
\begin{align*}
T_\boldk(\boldz) \coloneqq \prod_{j=1}^{d} T_{k_j}(z_j) \quad\text{with}\quad T_{k_j}(z_j)= \begin{cases} 1 & k_j = 0 \\
 \sqrt{2} \cos(k_j \arccos(z_j)) & k_j \not= 0 \end{cases} 
\end{align*}
being multivariate Chebyshev polynomials of dimension $d = N+1$, unknown coefficients $u_\boldk \in \C$ and an unknown, sparse index set $I \subset \N^{N+1}$. Note that the sparsity is a crucial requirement here, since a high-dimensional approximation of such form for non-sparse index sets $I$ is almost always computationally unfeasible due to the curse of dimension. On the other hand, sparsity occurs naturally based on the way we discretize the function $f$ as we show below. The behavior, that the solution is dominated by main effects and some low-order interactions, will reappear in our numerical experiments as a result of our general approach. A similar effect is known as the sparsity-of-effects principle in other fields of analysis, cf. \cite{WuHa21,HaSchShToTrWa22}.

Due to the simple structure of $\L$ and the initial condition, we now know that
\begin{align*}
u(x,\boldf) = f_0 x + \sum_{\ell=1}^{N-1} \frac{f_\ell} {2\pi\ii \ell} \e^{2\pi\ii \ell x},
\end{align*}
which can be rewritten using the univariate Chebyshev polynomials $T_{0}(z) = 1$ and $T_{1}(z) = \sqrt2 z$ as
\begin{align*}
u(x,\boldf) = \frac12 T_1(x) T_1(f_0) \prod_{j=1}^{N-1} T_0(f_j) + \sum_{\ell=1}^{N-1} \e^{2\pi\ii \ell x} \frac{T_1(f_\ell)} {2 \sqrt2 \pi\ii \ell} \prod_{\substack{j=0\\j\not=\ell}}^{N-1} T_0(f_j),
\end{align*}
since all the $f_\ell$ appear only linearly and decoupled. From this formula, we can directly read the structure of the index set $I$ as
\begin{align*}
I = \left\lbrace \begin{bmatrix}
1 \\ 1 \\ 0 \\ \vdots \\ 0
\end{bmatrix}, \begin{bmatrix}
0 \\ 0 \\ 1 \\ \vdots \\ 0
\end{bmatrix}, \ldots, \begin{bmatrix}
0 \\ 0 \\ 0 \\ \vdots \\ 1
\end{bmatrix}, \begin{bmatrix}
1 \\ 0 \\ 1 \\ \vdots \\ 0
\end{bmatrix}, \ldots, \begin{bmatrix}
1 \\ 0 \\ 0 \\ \vdots \\ 1
\end{bmatrix}, \begin{bmatrix}
2 \\ 0 \\ 1 \\ \vdots \\ 0
\end{bmatrix}, \ldots, \begin{bmatrix}
2 \\ 0 \\ 0 \\ \vdots \\ 1
\end{bmatrix}, \begin{bmatrix}
3 \\ 0 \\ 1 \\ \vdots \\ 0
\end{bmatrix}, \ldots \right\rbrace
\end{align*}
with the first row corresponding to the spatial dimension $x$ and the remaining rows corresponding to the coefficient dimensions $f_0,\ldots,f_9$ (in this order).

Therefore, we in fact know the index set $I$ and can compute the corresponding basis coefficients $u_\boldk$ simply by using e.g. Chebyshev rank-1 lattices \cite{Kae25}. So if the structure of a suitable index set $I$ in \eqref{eq:intro_sol} can be computed exactly, the hardest part of the approximation problem is already done and allows us to use high-dimensional cubature methods like Monte Carlo (MC) or Quasi-Monte Carlo (QMC) methods to derive the coefficients $u_\boldk$. A similar approach in the Fourier setting for a more complicated differential equation was investigated in \cite{GrIw24} and showed great results, even for very high dimensions. Further, the parametrization of the right-hand side function $f$ as well as the representation of the solution $u$ by its basis coefficients $u_\boldk$ is closely related to the concept of so-called encoder and decoder mappings as for example used in \cite{HeSchwaZe24}. There the existence of reasonable neural and spectral operators wrapped by such encoder and decoder mappings as well as their error bounds were studied. The fundamental difference of our approach is the missing decoder, since we use the Fourier coefficients $\boldf$ as the encoding of the function $f$ directly to compute the solution $u$ via \eqref{eq:intro_sol}.

Unfortunately, a direct analysis of the structure of the a priori unknown index set $I$ is often extremely difficult or simply impossible due to the complexity of the considered differential problem. Thus, we use a dimension-incremental algorithm presented in \cite{KaPoTa23} instead, using point samples to detect a suitable index set $I$ adaptively. More precisely, our general aim is to solve the high-dimensional approximation problem of approximating \eqref{eq:intro_sol} by detecting a reasonable and sparse index set $I$ containing the indices $\boldk$ corresponding to the largest (in absolute value) coefficients $u_\boldk$. Thus, we are using samples $u(x^{(j)},\boldf^{(j)})$ of the solution $u$, which we obtain by solving the differential equation \eqref{eq:intro_pde} for the fixed parameters $\boldf^{(j)}$. The adaptive detection of a good index set $I$ is made possible by the usage of adaptively chosen sampling points $(x^{(j)},\boldf^{(j)})$ and hence training data in the algorithm, which is one of the biggest differences of this approach to several deep learning techniques for operator learning associated with PDEs, where random training data is assumed.

The dimension-incremental algorithm in \cite{KaPoTa23}, based on the method from \cite{KaPoVo17}, works in arbitrary bounded orthonormal product bases (BOPB) and is capable of computing proper approximations with satisfying error bounds based on the underlying cubature method as shown for the Fourier case in \cite{BaTa23}. To generate the necessary point samples of the solution $u$, we will utilize classical differential equation solvers like the FEM. Once we detected the index set $I$, we can analyze its structure to gain insight on the interactions between the different parameters and variables and their influence on the solution $u$. Further, we can generalize the structure of the index set $I$ to proceed to even higher-dimensional or more detailed versions of the differential problem with a reasonable a priori guess, which coefficients $u_\boldk$ will be important therein.

This leads to an efficient and interpretable approximation of the solution $u$, expressed similarly to \eqref{eq:intro_sol}, allowing direct evaluations of $u$ for any right-hand side $f$ that is either of the same form or can be well approximated by it—without repeatedly solving the PDE. This operator-based perspective on PDE solvers offers an appealing alternative to classical model-based approaches, in which the structure of the differential operator is explicitly exploited to construct efficient solution algorithms. Our method can be seen as a synthesis of model-based and data-driven approaches: although we rely on sample-based recovery, the use of a structured index set and the decomposition into interpretable basis functions connects our framework closely to analytical PDE solvers, where solution structures are derived symbolically or semi-analytically. As such, it provides both interpretability and generalizability, thereby addressing one of the major limitations of many purely machine learning-based methods. In contrast to neural operator approaches or physics-informed neural networks (PINNs), which often operate as black-box models with limited insight into the underlying solution structures, our method maintains a clear connection to the analytical properties of the PDE and offers explicit control over the approximation process through the choice of basis and index set.

The remainder of this paper is organized as follows: In Section \ref{sec:theory}, we properly introduce our notations and the theoretical framework for the application of the dimension-incremental algorithm from \cite{KaPoTa23}, which is also explained briefly. Section \ref{sec:numerics} then investigates several test examples to present the application of our method. These examples include the Poisson and heat equation, where we have access to an analytical solution for comparison of the results as well more advanced differential problems like parametric diffusion equations with random coefficients and the non-linear Burgers' equation. Finally, we give a brief conclusion in Section \ref{sec:conclusion}.

The Python code is available at \url{https://github.com/fabiantaubert/nabopb}. It contains the dimension-incremental algorithm, the necessary rank-1 lattice cubature methods and the differential equation applications.

\section{Theory}
\label{sec:theory}

\subsection{Setting}

We investigate differential equations of the general form
\begin{align} \label{eq:general}
\L u = f
\end{align}
with a differential operator $\L:\U\rightarrow\F$, a right-hand side $f\in \F$ and the corresponding solution $u \in \U$. Here, $\U$ and $\F$ are suitable function spaces defined on the Lipschitz smooth $d$-dimensional spatial domain $\Omega \subset \R^d$. A typical example for the spaces $\U$ and $\F$, which will also appear in our numerical experiments in Section \ref{sec:numerics}, are the Sobolev space $H^1(\Omega)$ and its dual $H^{-1}(\Omega)$. Furthermore, together with boundary conditions like Dirichlet or Neumann conditions, the solution $u$ can be defined even on the closure $\overline{\Omega}$. 

To ensure that the differential problem \eqref{eq:general} is well-posed, we require certain conditions for the existence and uniqueness of the solution. A problem is said to be well-posed if there is a solution, the solution is unique, and the solution depends continuously on the data. These conditions are typically guaranteed through appropriate assumptions on the operator $\L$, the function spaces $\U$ and $\F$, and the boundary conditions.

For example, in the case where $\U=H^1(\Omega)$ and $\F=H^{-1}(\Omega)$, the Lax-Milgram theorem can be applied for elliptic problems to guarantee the existence and uniqueness of the solution under suitable conditions on the bilinear form associated with $\L$. In other cases, the Fredholm alternative or other functional analysis results may be used to establish these properties. Additionally, the regularity of the solution often depends on the smoothness of the data and the boundary conditions, which in turn influences the choice of suitable function spaces for approximation.

\subsection{Sparse Approximation}

Assuming that \eqref{eq:general} has a unique solution for each right-hand side function $f \in \F$, our goal is to learn the solution mapping $\G: \F \rightarrow \U$, ending up with an apporoximation $\tilde{\G}$ of $\G$.

In order to do so we start by approximating the function $f$ by 
\begin{align} \label{eq:f_approx}
f(\boldx) \approx \sum_{j=1}^n a_j A_j(\boldx) & & \boldx \in \Omega,
\end{align}
with some fixed functions $A_j, j = 1,\ldots,n,$ and coefficients $\bolda = (a_1,\ldots,a_n)\in\C^n$. Favorable approximative parametrizations of the function $f$ by the coefficients $\bolda$ should be accurate and efficient, i.e., the possible error should be reasonably bounded and there exist fast methods to compute the coefficients $\bolda$ for a given function $f$ and vice versa.

Further, we assume the solution $u: \D \rightarrow \C$ to have a basis expansion of the form
\begin{align} \label{eq:basexp_of_u}
u(\boldx,\bolda) \coloneqq \sum_{\boldk \in \N^{d+n}} c_{\boldk} \Phi_\boldk(\boldx,\bolda) & & (\boldx,\bolda) \in \D,
\end{align}
with $\lbrace \Phi_\boldk: \boldk \in \N^{d+n}\rbrace$ a bounded orthonormal product basis (BOPB) in the separable Hilbert space $\H = L_2(\D,\mu)$ on the Cartesian product type domain $\D$ and basis coefficients $c_\boldk \in \C, \boldk\in\N^{d+n}$. See \cite[Sec.~1.1]{KaPoTa23} for more details on the notion of a BOPB and the corresponding domains and spaces. Further, see Remark \ref{rem:transform} for the relation between the domains $\D$ and $\overline{\Omega} \times C^n$, which would be the canonical domain for $u$ since $\boldx \in \overline{\Omega}$ and $\bolda \in \C^n$.

Having \eqref{eq:basexp_of_u}, we aim to approximate $u$ by a truncation
\begin{align} \label{eq:basexp_of_u_trunc1}
S_I u(\boldx,\bolda) \coloneqq \sum_{\boldk \in I} c_{\boldk} \Phi_\boldk(\boldx,\bolda),
\end{align}
with some a priori unknown index set $I \subset \N^{d+n}, |I| < \infty$, followed by an approximation
\begin{align} \label{eq:basexp_of_u_trunc}
S_I^\A u(\boldx,\bolda) \coloneqq \sum_{\boldk \in I} \hat{u}_{\boldk} \Phi_\boldk(\boldx,\bolda),
\end{align}
where $\hat{u}_{\boldk} \in \C, \boldk \in I,$ are approximations of the true coefficients $c_{\boldk}$. Note that the detection of a ``good" index set $I$ in general leads to a non-linear approximation problem. With \eqref{eq:basexp_of_u_trunc} we then have an approximation of the solution mapping: For every right-hand side $f$ we determine the coefficients $\bolda$ and plug them into $S_I^\A u$, which yields us an explicit representation of an approximation of $u$. So in fact, combining the discretization mapping $f \mapsto \bolda$ with the approximation mapping $\bolda \mapsto S_I^\A u(\cdot,\bolda)$ yields us an approximation $\tilde{G}$ of the solution map $G: \F \rightarrow \U$. Furthermore, the structure of the index set $I$ as well as the size of the corresponding approximated coefficients $\hat{u}_\boldk$ may reveal interesting insights on the structure of the solution $u$ and its dependence on the coefficients $\bolda$ and, equivalently, the dependence on the right-hand side function $f$.

Estimating the error of the approximation in the $L_\infty$ norm and using the boundedness of our basis functions $\Phi_\boldk$, we get:
\begin{align*}
\norm{u-S_I^\A u}_\infty &\leq \norm{u - S_I u}_\infty + \norm{S_I u - S_I^\A u}_\infty \\
&= \norm{\sum_{\boldk \not\in I} c_{\boldk} \Phi_\boldk}_\infty + \norm{\sum_{\boldk \in I} (c_{\boldk}-\hat{u}_\boldk) \Phi_\boldk}_\infty \\
&\leq \sum_{\boldk \not\in I} \abs{c_{\boldk}} \norm{\Phi_\boldk}_\infty + \sum_{\boldk \in I} \abs{c_{\boldk}-\hat{u}_\boldk} \norm{ \Phi_\boldk}_\infty \\
&\leq B \left( \sum_{\boldk \not\in I} \abs{c_{\boldk}} + \sum_{\boldk \in I} \abs{c_{\boldk}-\hat{u}_\boldk} \right)
\end{align*}
Note that the boundedness constant $B$ depends on the BOPB and therefore on the space $\H$ of our approximation. As an example, consider the well-known Fourier system using trigonometric polynomials $\exp(2\pi\ii\langle\boldk,\cdot\rangle)$ as basis functions, which comes with the boundedness constant $B=1$ independent of the dimenson $d$. For the Chebyshev basis already introduced in Section~\ref{sec:introduction}, the boundedness constant is $B=2^{d/2}$. If the approximation problem considered has a limited number of interacting dimensions $d_s \leq d$, the effective boundedness constant can be replaced with $B=2^{d_s/2}$ instead.

The terms inside the brackets, which we will refer to as the truncation error $\sum_{\boldk \not\in I} \abs{c_{\boldk}}$ and the coefficient approximation error $\sum_{\boldk \in I} \abs{c_{\boldk}-\hat{u}_\boldk}$, are mainly influenced by the index set $I$ and the approximated coefficients $\hat{u}_\boldk$. Hence, we need not only to compute good approximations of the coefficients $c_{\boldk}$ but also to detect a ``good" index set $I$, hopefully containing the largest (in terms of absolute value) coefficients $c_{\boldk}$, to make \eqref{eq:basexp_of_u_trunc} a reasonable approximation of the solution $u$. 

\subsection{The dimension-incremental method}

In the present work, we will use the nonlinear approximation method for high-dimensional function approximation proposed in \cite{KaPoTa23} in order to receive the desired approximation \eqref{eq:basexp_of_u_trunc}. First, we summarize some main aspects of the dimension-incremental algorithm here and refer to \cite[Sec.~2]{KaPoTa23} for more detailed explanations and proper definitions of the used notations. A simplified version of the algorithm is also given in Algorithm \ref{alg:main}.

Suppose for now for simplicity that we are interested in the approximation of a $d$-dimensional target function $g: \D \rightarrow \C$ of the form $g = \sum_{\boldk\in I} \ghat_\boldk\Phi_\boldk$ with the unknown index set $I \subset \N^d$. Later, the target function $g$ will be the solution $u$ in $d+n$ dimensions. Motivated by the estimate above, we aim for an $s$-sparse index set $I$, i.e., we have $|I| = s$, corresponding to basis coefficients $c_\boldk$ with large absolute values. 

Roughly spoken, the algorithm uses samples of $g$ to detect reasonable indices $k_j$ of $\boldk = (k_j)_{j=1}^{d}, \boldk \in I,$ in each dimension $j=1,\ldots,d$ and reasonable combinations thereof. In order to do so, only a search space $\Gamma \supset I$ is needed in advance. Commonly, we choose search spaces like the full (non-negative) grid $\Gamma = [0,N]^d$ with some parameter $N \in \N$, which we will call extension from now on. If there is additional initial knowledge on the structure of the desired index set $I$, the choice of $\Gamma$ can be improved.

The algorithm starts by investigating the one-dimensional projections $\mathcal{P}_{\{t\}}(\Gamma) \coloneqq \{k \in \N\,|\,\exists \boldk \in \Gamma: \boldk_t=k\}$ for all $t=1,\ldots,d$ by constructing a suitable cubature rule for integrals of the form 
\begin{align}\label{eq:proj_coef_1d}
\ghat_{\{t\},k_t}(\tilde{\boldx}) \coloneqq \int g(\xi, \tilde{\boldx})_{\{t\}} \overline{\Phi_{\{t\},k_t}(\xi)} \diff \xi
\end{align}
with $\Phi_{\{t\},k_t}$ the one-dimensional basis function of the $t$-th dimension of our BOPB. The notation $g(\xi, \tilde{\boldx})_{\{t\}}$ refers to sampling  values of $g$ using $\xi$ in the $t$-th dimension and $\tilde{\boldx}$ for the remaining dimensions. The algorithm then computes these so-called projected coefficients $\ghat_{\{t\},k_t}$ using this cubature rule and samples of the target function $g$ for a particular random anchor $\tilde{\boldx}$. The absolute value of these projected coefficients $\ghat_{\{t\},k_t}$ can be seen as an indicator whether or not $k_t$ is important, i.e., if $k_t$ should appear in the $t$-th component of any index $\boldk \in I$. Hence, the algorithm takes the sparsity $s$ largest projected coefficients fulfilling $|\ghat_{\{t\},k_t}| \geq \delta_+$ for some initially chosen detection threshold $\delta_+$ and adds the corresponding indices $k_t$ to a temporary index set $I_{\{t\}}$. Since the computation of the projected coefficients $\ghat_{\{t\},k_t}$ involves randomness due to the randomly drawn anchor $\tilde{\boldx}$, this computation is repeated $r$ times with $r$ being the number of detection iterations with different anchors $\tilde{\boldx}^{(j)}, j=1,\ldots,r$.

The temporary index sets $I_{\{t\}}$ with the reasonable indices for each dimension $t=1,\ldots,d$ are now combined to proceed in a dimension-incremental way. Starting with $t=2$ a new candidate set $K:=(I_{\{1,\ldots,t-1\}} \times I_{\{t\}})\cap\mathcal{P}_{\{1,\ldots,t\}}(\Gamma)$ is formed, containing now higher-dimensional indices $\boldk \in \N^{|\{1,\ldots,t\}|}$. Again, a suitable $t$-dimensional cubature method, e.g., multiple rank-1 lattices as in \cite{Kae16}, is constructed and evaluated using samples of the target function $g$ to compute the projected coefficients
\begin{align*}
\ghat_{\{1,\ldots,t\},\boldk}(\tilde{\boldx}) \coloneqq \int g(\boldxi, \tilde{\boldx}) \overline{\Phi_{\{1,\ldots,t\},\boldk}(\boldxi)} \diff \boldxi
\end{align*}
for the indices $\boldk \in K$, which is the natural generalization of \eqref{eq:proj_coef_1d} to multiple dimensions $\{1,\ldots,t\}$. As before, the algorithm collects those indicies $\boldk \in \N^t$ with the sparsity $s$ largest (in absolute value) projected coefficients $\ghat_{\{1,\ldots,t\},\boldk}$ in the temporary index set $I_{\{1,\ldots,t\}}$. These indices are now the tuples, which can still appear in the first $t$ components of the indices in the final index set $I$. If $t<d$ holds, this process is again influenced by the randomly chosen anchor $\tilde{\boldx}$ and therefore repeated $r$ times. This process is then repeated with $t+1$ instead of $t$ until $t=d$, where the projected coefficients $\ghat_{\{1,\ldots,d\},\boldk}$ do not longer depend on a random anchor $\tilde{\boldx}$ at all. We finally set $I=I_{\{1,\ldots,d\}}$ and $\ghat_\boldk \coloneqq \ghat_{\{1,\ldots,d\},\boldk}$ for all $\boldk\in I_{\{1,\ldots,d\}}$. The final output of the algorithm is the desired index set $I$ as well as approximations $\ghat_\boldk$ of the true basis coefficients for each $\boldk \in I$. 

A crucial requirement for this algorithm is the possibility to access the sampling values $g(\boldx)$ for arbitrary $\boldx$ during the algorithm, e.g. by a black-box function handle of $g$. This is since the necessary sampling points $\boldx$, for which the corresponding sampling values $g(\boldx)$ are needed, are not known a priori, but are computed adaptively during the algorithm based on the current candidate sets $K$ and the constructed cubature methods, combined with the random anchors $\tilde{\boldx}$. We again encourage the reader to see \cite{KaPoTa23} for a more detailed and rigorous explanation of this concept, the theorem on the theoretical detection guarantee and simple numerical examples as well as several comments and discussions on the capabilities and restrictions of this algorithm. Further, similar approaches as in \cite[Sec.~3.1.2]{zech18} for the selection of the largest projected coefficients $\ghat_{\{1,\ldots,t\},\boldk}(\tilde{\boldx})$ of the candidate sets $K$ in each step might speed up these detections and therefore the whole algorithm, if additional information on the behavior of these projected coefficients is known.

\begin{algorithm}[t]
	\caption{Dimension-incremental Algorithm (Simplified)}\label{alg:main}
  \begin{small}
	\begin{tabular}{p{1.3cm}p{2.4cm}p{10.6cm}}
		Input:  %
		& $\Gamma\subset\N^d$ \hfill & search space \\
		& $g$ & target function $g$ as black box (function handle) \\
		& $s\in\N$ & sparsity parameter \\  
		& $\delta_+ > 0$ & detection threshold \\
		& $r\in\N$ & number of detection iterations%
	\end{tabular}
	\begin{algorithmic}
		\item[(Step 1)] [Single component identification]
		\STATE {\bfseries for} $t:=1,\ldots,d$ {\bfseries do}
		\STATE \hspace{1em} Set $I_{\{t\}}:=\emptyset$.
		\STATE \hspace{1em} Compute a suitable cubature method for $\mathcal{P}_{\{t\}}(\Gamma)$.
		\STATE \hspace{1em} {\bfseries for} $i:=1,\ldots,r$ {\bfseries do}
		\STATE \hspace{2em} Draw a random anchor $\tilde{\boldx}$.
		\STATE \hspace{2em} Sample $g$ at the necessary sampling points (the cubature nodes combined with $\tilde{\boldx}$).
		\STATE \hspace{2em} Compute the projected coefficients $\ghat_{\{t\},k_t}(\tilde{\boldx})$ for $k_t \in \mathcal{P}_{\{t\}}(\Gamma)$.
  		\STATE \hspace{2em} Add the (up to) $s$ indices $k_t$ with the largest proj. coef. $|\ghat_{\{t\},k_t}(\tilde{\boldx})| \geq \delta_+$ to $I_{\{t\}}$.
		\STATE \hspace{1em} {\bfseries end for} $i$
		\STATE {\bfseries end for} $t$
		\item[(Step 2)] [Coupled component identification]
		\STATE {\bfseries for} $t:=2,\ldots,d$ {\bfseries do}
		\STATE \hspace{1em} If $t<d$, set $\tilde r:=r$ and otherwise $\tilde r:=1$. 
		\STATE \hspace{1em} Set $I_{\{1,\ldots,t\}}:=\emptyset$.
		\STATE \hspace{1em} Construct the index set $K:=(I_{\{1,\ldots,t-1\}} \times I_{\{t\}})\cap\mathcal{P}_{\{1,\ldots,t\}}(\Gamma)$.
		\STATE \hspace{1em} Compute a suitable cubature method for $K$.
		\STATE \hspace{1em} {\bfseries for} $i:=1,\ldots,\tilde r$ {\bfseries do}
		\STATE \hspace{2em} Draw a random anchor $\tilde{\boldx}$.
		\STATE \hspace{2em} Sample $g$ at the necessary sampling points (the cubature nodes combined with $\tilde{\boldx}$ if $t<d$).
		\STATE \hspace{2em} Compute the projected coefficients $\ghat_{\{1,\ldots,t\},\boldk}(\tilde{\boldx})$ for ${\boldk \in K}$.
		\STATE \hspace{2em} Add the (up to) $s$ indices $\boldk$ with the largest proj. coef. $|\ghat_{\{1,\ldots,t\},\boldk}(\tilde{\boldx})| \geq \delta_+$ to $I_{\{1,\ldots,t\}}$.
		\STATE \hspace{1em} {\bfseries end for} $i$
		\STATE {\bfseries end for} $t$
		\item[(Step 3)]
		\STATE Set $I:=I_{\{1,\ldots,d\}}$ and $\ghat_\boldk \coloneqq \ghat_{\{1,\ldots,d\},\boldk}$ for all $\boldk\in I_{\{1,\ldots,d\}}$.
	\end{algorithmic}

	\begin{tabular}{p{1.3cm}p{2.4cm}p{10.6cm}}
		Output: & $I\subset\Gamma\subset\N^d$ & detected index set\\
		& $(\ghat_\boldk)_{\boldk\in I} \in\C^{\abs{I}}$ & approximated coefficients with $|\ghat_{\boldk}| \geq \delta_+$ \\
	\end{tabular}
  \end{small}
\end{algorithm}

\subsection{Black-box sampling the differential equation}
\label{subsec:theory_bbs}

For our application of Algorithm \ref{alg:main}, we have $u$ as the target function and hence sampling points of the form $(\boldx,\bolda)$. Each ``sample" $u(\boldx,\bolda)$ is then the value of the solution $u$ of \eqref{eq:general} for a given parameter $\bolda$, evaluated at the spatial point $\boldx$. In order do receive such samples, we utilize numerical solvers, e.g. the finite element method. With the given parameter $\bolda$, we can compute the (fixed) right-hand side function $f$, solve the differential equation for this particular $f$ numerically and evaluate the approximated solution at $\boldx$.

Since the algorithm only requests the final sampling value $u(\boldx,\bolda)$, the choice of the particular method or numerical solver for the differential equation is completely free, as long as the approximated sampling value we compute and return to the algorithm is a reasonable approximation of the true value $u(\boldx,\bolda)$. This non-intrusive behavior of our algorithm is the reason for its generality. The properties of the differential equation \eqref{eq:general} as well as possible difficulties therein are mostly dealt with by the numerical solver. As long as there is any possibility to obtain reasonable estimates of $u(\boldx,\bolda)$ for a given $(\boldx,\bolda)$, we can plug this method into our black-box sampling step and are done. Note that while the accuracy and efficiency of the numerical solver used obviously directly affect the accuracy and efficiency of our dimension-incremental method, we will not investigate the properties of these solvers in more detail in this work.

\begin{remark} \label{rem:transform}
Generally, the product type domain $\D$, where we can apply Algorithm \ref{alg:main}, will not coincide with the domain $\overline{\Omega} \times \C^n$ of our solution $u$. Hence, we need to transform and or restrict this domain carefully, such that the sampling points given by our algorithm are suitable for the differential operator.

Since $\overline{\Omega}$ will be some compact domain for many applications, it is often enough to apply a simple transformation $\mathcal{T}$ for the spatial variable $\boldx$, e.g., the continuous and bijective linear transformation $\mathcal{T} \boldx = m_1 \boldx + m_2 \boldone$ with two constants $m_1, m_2 \in \R$, mapping $\boldx$ to the desired domain $\overline{\Omega}$.

On the other hand, the parameters $\bolda \in \C^n$ are more difficult to handle. In this case, we will often have to restrict the domain of $\bolda$ to e.g. some compact interval again, before thinking about possible transformations as we did for the spatial part. This is obviously a loss of generality and the restriction needs to be performed carefully, such that most reasonable source functions $f$ can still be approximated well enough using the restricted $\bolda$.

For examples of such transformations and restrictions, we refer to the particular examples in the following section.
\end{remark}

\begin{remark}\label{rem:parameters}
While we stick to the mentioned setting \eqref{eq:general} for the theoretical part of this paper to preserve clarity in the notations, our numerical experiments in Section \ref{sec:numerics} also include some further variations, which we will only briefly mention in this remark.

First, we can also consider parametrized differential operators $\L_\boldtheta$ with some parameter $\boldtheta \in \R^{n_\theta},n_\theta\in\N,$ and the corresponding solution mapping $\G: \F \times \R^n \rightarrow \U$. Correspondingly, \eqref{eq:basexp_of_u} then becomes
\begin{align}\label{eq:param_u_approx}
u(\boldx,\bolda,\boldtheta) \coloneqq \sum_{\boldk \in \N^{d+n+n_\theta}} c_{\boldk} \Phi_\boldk(\boldx,\bolda,\boldtheta)
\end{align}
with another separable Hilbert space $\H$ and corresponding BOPB $\lbrace \Phi_\boldk: \boldk \in \N^{d+n+n_\theta}\rbrace$. The truncated and approximated version \eqref{eq:basexp_of_u_trunc} is then modified in the same way, now with an unknown index set $I \subset \N^{d+n+n_\theta}$. An example of this variation can be found in Section \ref{subsec:rand_diff_eq}.

Similarly, we can consider time-dependent differential operators $\L$ with respect to some time variable $\tau\in [0,T]$ and their corresponding solution mapping $\G: \F \times [0,T] \rightarrow \U$. As before, we end up with the representation
\begin{align*}
u(\boldx,\tau,\bolda) \coloneqq \sum_{\boldk \in \N^{d+1+n_\theta}} c_{\boldk} \Phi_\boldk(\boldx,\tau,\bolda)
\end{align*}
and proceed similarly as above, now with the $d+1+n$-dimensional separable Hilbert space $\H$ and an unknown index set $\I \subset \N^{d+1+n}$. 

In each case, we proceed to the approximation $S_I^\A u$ from \eqref{eq:basexp_of_u_trunc}. This time, the analysis of the index set $I$ and the coefficients $\hat{u}_\boldk$ can give additional information about the dependence and interaction of the spatial variable $\boldx$ and the right-hand side $f$ not only with each other, but also with the parameters $\boldtheta$ or the time variable $t$. 

Obviously, a combination of these two variations, i.e.\ a parameter- and time-dependent differential equation, can be treated in an analogous way, cf.\ Section \ref{subsec:heat_eq} with the one-dimensional heat equation.
\end{remark}

\section{Numerics}
\label{sec:numerics}

In this Section, we will test our approach on several test problems, such as multiple diffusion equations and Burgers' equation, and discuss the results. We show, that our approach leads to sparse index sets $I$ that can be used directly for the high-dimensional approximation of the solution $u$ or further generalized to even higher-dimensional problems. We will briefly investigate such a generalization for the first model example at the end of Section \ref{subsubsec:poisson_1d_fc}. For the other examples, we focus on the detection of a suitable index set $I$ and omit the generalization to higher dimensions, since the first part is the main goal of Algorithm \ref{alg:main}.

As mentioned in Section \ref{sec:theory}, we need to choose a suitable BOPB for the solution $u$ to achieve the basis expansion \eqref{eq:basexp_of_u}. In all of the following examples, we will work with the tensorized Chebyshev polynomials
\begin{align*}
T_\boldk(\boldz) \coloneqq \prod_{j=1}^{d+n} T_{k_j}(z_j) \quad\text{with}\quad T_{k_j}(z_j)= \begin{cases} 1 & k_j = 0 \\
 \sqrt{2} \cos(k_j \arccos(z_j)) & k_j \not= 0 \end{cases} 
\end{align*}
on the domain $\D = [-1,1]^{d+n}$ as used in \cite{KaPoTa23}. Hence, the approximation \eqref{eq:basexp_of_u_trunc} inserting $\boldz \coloneqq (\boldx,\bolda)$ with $\boldx \in \R^d$ and $\bolda \in \C^n$ becomes
\begin{align} \label{eq:cheby_approx_of_u}
S_I^\A u(\boldx,\bolda) = \sum_{\boldk \in I} \hat{u}_{\boldk} T_\boldk(\boldx,\bolda),
\end{align}
where $I \subset \N^{d+n}$ is the detected index set and $\hat{u}_{\boldk}$ are the approximations of the corresponding exact basis coefficients $c_\boldk$ from \eqref{eq:basexp_of_u}.

To investigate the accuracy of our method, we consider for a fixed coefficient $\bolda \in \C^n$ the relative $\ell_2$-error
\begin{align} \label{eq:err}
\text{err}(\bolda) \coloneqq \frac{\norm{S_I^\A u(\boldx,\bolda) - u(\boldx,\bolda)}_{\ell_2}}{\norm{u(\boldx,\bolda)}_{\ell_2}} = \frac{\left(\sum_{j=1}^{G}\abs{S_I^\A u(\boldx^{(j)},\bolda) - u(\boldx^{(j)},\bolda)}^2\right)^{\frac12}}{\left(\sum_{j=1}^{G}\abs{u(\boldx^{(j)},\bolda)}^2\right)^{\frac12}},
\end{align}
where $\boldx^{(j)}$ for $j=1,\ldots,G$ are equidistant grid points in the spatial domain $\Omega$. We then proceed by computing this error for numerous, randomly drawn coefficients $\bolda$ and investigating the corresponding range as well as the first quartile, the median and the second quartile of this statistical test (dividing the results into four equal parts).

All tests are performed in Python\textsuperscript{TM} and can be found together with the algorithm in \cite{gitNABOPB}. The overall runtime depends significantly on the specific problem, the desired accuracy, the chosen parameters, and the available computational resources. Therefore, the runtime for each example is briefly discussed in the corresponding sections. To accelerate the sampling process, we used parallelization with up to $180$ separate workers, but not GPU computing. If not stated otherwise, the dimension-incremental algorithm uses the following parameters and settings: \begin{itemize}
\item the cubature method: Chebyshev multiple rank-1 lattices as described in \cite[Sec.~4.2]{KaPoTa23}
\item the search space $\Gamma$: (non-negative) full grid $[0,N]^{d+n}$ in $d+n$ dimensions with extension $N$ and no superposition assumption
\item the detection threshold $\delta_+ = 10^{-12}$
\item the number of detection iterations $r=5$.
\end{itemize}
The sparsity $s$ will be given for each test separately. See \cite{KaPoTa23} for more detailed information on these parameters and settings and how they affect the behavior of the algorithm.

\subsection{The one-dimensional Poisson equation}
\label{subsec:poisson_1d}

The following example considers a rather simple differential equation in order to demonstrate the application of our proposed method for the first time.

Given a source function $f: (0,1) \rightarrow \C$, the one-dimensional Poisson equation with homogeneous Dirichlet boundary conditions reads as
\begin{equation} \label{eq:poisson_1d}
\begin{aligned} 
-\frac{d^2}{d x^2} u(x) &= f(x), & x \in (0,1), \\
u(0) = u(1) &= 0. &
\end{aligned}
\end{equation}
For this differential operator $\L = -\frac{d^2}{d x^2}$ and these particular boundary conditions, we go with the usual choice of function spaces $\U = H_0^1((0,1))$ and $\F = H^{-1}((0,1))$, the dual of $H_0^1((0,1))$. As described in Section \ref{sec:theory}, our first step is to find a suitable parametrization \eqref{eq:f_approx} of the function $f$. For this first example, we will consider two different approaches here: \begin{itemize}
\item a parametrization of the source function $f$ by its first Fourier coefficients,
\item a parametrization of the source function $f$ by a B-Spline approximation.
\end{itemize}
In all of these cases we restrict ourselves to a discretization of $f$ using only $n$ parameters. Together with the spatial dimension $d=1$ this results in a $n+1$-dimensional approximation problem. While choosing a larger $n$ should lead to more accurate approximations of $f$ and thus to an overall better quality of the approximation $S_I^\A u$ for general $f$, the additional dimensions will also result in higher sampling and computational complexities. Therefore, we have to choose reasonable limits for $n$ in the upcoming examples.

\subsubsection{Fourier series parametrization} \label{subsubsec:poisson_1d_fc}
We consider a parametrization of the source function $f$ by its first $n \in 2\N+1$ Fourier coefficients $\bolda = (a_{-\frac{n-1}{2}},\ldots,a_{\frac{n-1}{2}}) \in \C^n$, i.e.,
\begin{align}\label{eq:pois1d_f_four}
f(x) \approx \sum_{\ell=-\frac{n-1}{2}}^{\frac{n-1}{2}} a_\ell \e^{2\pi\ii \ell x}.
\end{align}
These Fourier coefficients $\bolda$ can be computed efficiently for reasonable functions $f$ using the well known FFT. Note that this approximation of $f$ will always be $1$-periodic, forcing the implicit assumption that the function $f$ is either a $1$-periodic function itself or can be well approximated by such a function up to some extend.

Using this truncated Fourier series as the right-hand side of the differential equation \eqref{eq:poisson_1d} the solution $u$ of the one-dimensional Poisson equation is then given analytically by
\begin{align} \label{eq:u_pois_anal}
u(x,\bolda) = \frac{a_0}{2}x(1-x) + \sum_{\substack{\ell=-\frac{n-1}{2} \\ \ell\neq 0}}^{\frac{n-1}{2}} \frac{a_\ell}{4\pi^2 \ell^2} (\e^{2\pi\ii\ell x}-1).
\end{align}
This formula can be used directly as the black-box sampling strategy necessary for our algorithm, see Section \ref{subsec:theory_bbs}, to generate the necessary sampling values $u(x^\ast,\bolda^\ast)$ for any sampling point $(x^\ast,\bolda^\ast)$. To demonstrate the general application of Algorithm \ref{alg:main}, we use this direct method to be able to neglect errors made when solving the differential equation numerically and focus on the approximation of the solution $u$ directly. 

As mentioned in Remark \ref{rem:transform}, we need to pay attention since the original domain $[0,1] \times \C^n$ doesn't match our function approximation domain $\D = [-1,1]^{n+1}$. For the spatial part, we apply the transformation $\mathcal{T} x=\frac12(x+1)$ to perform the shift between $[-1,1]$ and $[0,1]$. For simplicity, we directly assume the restriction $\bolda\in[-1,1]^n$ for the Fourier coefficients $\bolda$ such that we can omit further transformations. Note that this implies that we are only interested in right-hand side functions $f$, which can be well approximated by using such Fourier coefficients $\bolda$ during this artificial example. Overall, the final function, which we are going to approximate here, is now
\begin{align*}
\tilde{u}(\tilde{x},\bolda) \coloneqq u(\mathcal{T}\tilde{x},\bolda) = u\left(\frac12(\tilde{x}+1),\bolda\right) & & \tilde{x} \in [-1,1], \bolda \in [-1,1]^n.
\end{align*}
Using the explicit formula \eqref{eq:u_pois_anal}, we get
\begin{align} \label{eq:u_pois_anal_trans}
\tilde{u}(\tilde{x},\bolda) = \frac{a_0}{8}(1-\tilde{x}^2) + \sum_{\substack{\ell=-\frac{n-1}{2} \\ \ell\neq 0}}^{\frac{n-1}{2}} \frac{a_\ell}{4\pi^2 \ell^2} ((-1)^\ell \e^{\pi\ii\ell \tilde{x}}-1).
\end{align}
\begin{remark}
Note that the particular choice of the domain $\D$ and the corresponding basis of $\H$ are not unique. Since we are only restricting the Fourier coefficients $\bolda$ and not transforming them, the solution $u$ is obviously not periodic with respect to these variables, so our decision to use the tensorized Chebyshev polynomials for the approximation is reasonable here. However, we could have applied various transformations $\mathcal{T}$ to $\bolda$, including those that force a periodic dependence of $u$ on $\bolda$ such as the tent-transform, cf.~\cite{CoKuNuSu16}. Then, together with the periodicity in $x \in [0,1]$ due to the boundary conditions $u(0)=u(1)=0$, the solution $u$ would be periodic (but not smooth) in all $n+1$ dimensions. In such a scenario, we could use the high-dimensional torus domain $\D=\T^{n+1}$ as well as a Fourier basis for the approximation space $\H$.
\end{remark}
We use $n=9$ as the amount of Fourier coefficients $a_\ell$ for our tests, which results in the overall dimension $d+n=10$. Further, we choose the sparsity $s=1000$ and the extension $N=64$ of the search space $\Gamma$. Since we do not use a numerical solver but the exact solution, each solution sample can be generated in about $10^{-4}$ seconds, resulting in an overall runtime of our algorithm of about $2$ minutes.

\begin{figure}[t]
\centering
\begin{tikzpicture}
\begin{axis}[
	boxplot/draw direction=x,
	xmode=log,
	xmin=0.0000001,
	xmax=1,
	ytick={1,2,3,4,5},
	yticklabels={$[-1{,}1]$,$[-2{,}2]$,$[-4{,}4]$,$[-1{,}1] + [-1{,}1] \ii$, $[-2{,}2] + [-2{,}2] \ii$},
	xlabel={approximation error $\text{err}(\bolda)$},
	ylabel={range of the entries of $\bolda$},
	width=0.85\textwidth,
    height=0.40\textwidth,
	]
	\pgfplotstableread[col sep=comma]{1a_err.csv}\datatable;
	\addplot+[
  		boxplot prepared from table={
    			table=\datatable,
    			row=0,
    			lower whisker=lw,
    			upper whisker=uw,
    			lower quartile=lq,
    			upper quartile=uq,
    			median=med
  		}, my boxplot style,
    ] coordinates {};
    \addplot+[
  		boxplot prepared from table={
    			table=\datatable,
    			row=1,
    			lower whisker=lw,
    			upper whisker=uw,
    			lower quartile=lq,
    			upper quartile=uq,
    			median=med
  		}, my boxplot style,
    ] coordinates {};
    \addplot+[
  		boxplot prepared from table={
    			table=\datatable,
    			row=2,
    			lower whisker=lw,
    			upper whisker=uw,
    			lower quartile=lq,
    			upper quartile=uq,
    			median=med
  		}, my boxplot style,
    ] coordinates {};
    \addplot+[
  		boxplot prepared from table={
    			table=\datatable,
    			row=3,
    			lower whisker=lw,
    			upper whisker=uw,
    			lower quartile=lq,
    			upper quartile=uq,
    			median=med
  		}, my boxplot style,
    ] coordinates {};
    \addplot+[
  		boxplot prepared from table={
    			table=\datatable,
    			row=4,
    			lower whisker=lw,
    			upper whisker=uw,
    			lower quartile=lq,
    			upper quartile=uq,
    			median=med
  		}, my boxplot style,
    ] coordinates {};
    \end{axis}
\end{tikzpicture}
\caption{The relative approximation error $\text{err}(\bolda)$ for $10000$ randomly drawn $\bolda$ when using the Fourier series parametrization. The box-and-whisker plots show the median, the first and the second quartile as well as the maximal and minimal error observed. The five plots indicate different choices for the range of the Fourier coefficient $\bolda$, including two cases with complex-valued Fourier coefficients. The range $[-1,1]$ coincides with the training data used.}
\label{fig:pois1a_err}
\end{figure}
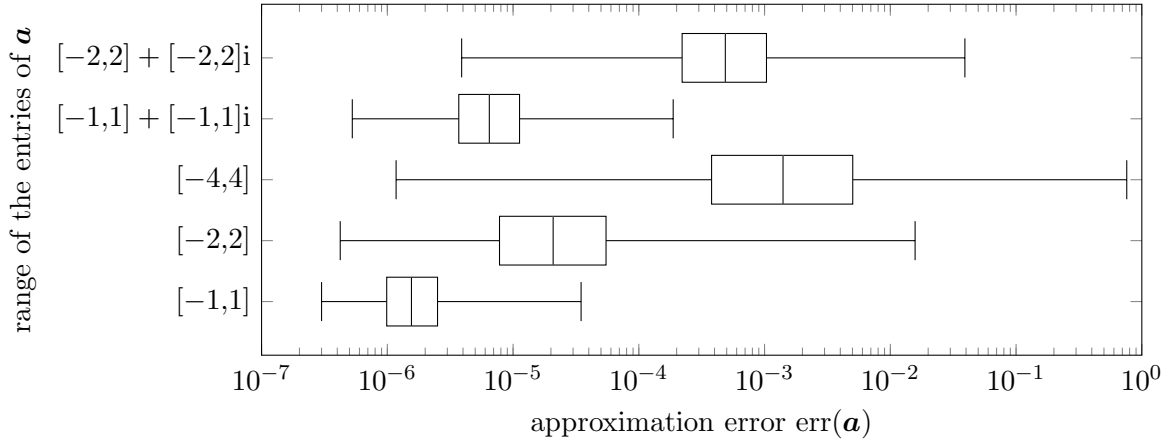

The accuracy of our approximation is shown in Figure \ref{fig:pois1a_err}. Therein, we used $10000$ randomly drawn Fourier coefficients $\bolda$ and computed the relative $\ell_2$-error $\text{err}(\bolda)$ using $G = 1000$ equidistant grid points in the spatial domain. We use box-and-whisker plots to illustrate the statistical distribution here, where the central line inside the box indicates the median. On each side of the median, the box contains 25\% of the data. Outside the box, the whiskers indicate the maximal and minimal error observed. Since we did not specify outliers in our data, the box-and-whisker plot truly covers the full range of observed errors.

The first plot with the range $[-1,1]$ is the true approximation error, since it used the same range for the entries of $\bolda$ as we used during our approximation. Although computed coefficients $\hat{u}_\boldk$ from \eqref{eq:cheby_approx_of_u} smaller than $10^{-7}$ are not necessarily true basis coefficients but mainly artifacts because of numerical errors, the overall approximation accuracy is still satisfactory. The other plots in Figure \ref{fig:pois1a_err} show results with larger or even complex domains for the test Fourier coefficients $\bolda$. For this transfer learning scenario it shows, that our approximation is also applicable for slightly larger domains of $\bolda$, and therefore more source functions $f$, than the restricted ones from the training setting. However, the further the test cases are from the training setting, the larger our relative approximation error $\text{err}(\bolda)$ grows. Thus, we strongly recommend using another restriction $\bolda \in [-\alpha,\alpha]^n$, $\alpha > 1$, and transforming them similarly to the spatial variable $x$ already in the training process, if desired functions $f$ require such Fourier coefficients $\bolda$.

\begin{figure}[t]
\centering
\begin{tikzpicture}
\pgfplotstableread{ind.dat}\mytable
\foreach \i in {0,...,39} {
    \draw[dashed, thin, gray, opacity=0.5] (\i*10 pt,0) -- (\i*10 pt,-90 pt);
}
\foreach \j in {0,...,9} {
    \draw[dashed, thin, gray, opacity=0.5] (0,-\j*10 pt) -- (390 pt,-\j*10 pt);
}
\foreach \i in {0,...,39}{
    \foreach \j in {0,...,9}{
    \pgfplotstablegetelem{\i}{\j}\of\mytable
    \ifnum\pgfplotsretval=0\relax\else
    \node[circle, minimum size=10pt, inner sep=0pt, fill=color6, opacity=0.\pgfplotsretval, text opacity= 0.75] at (\i*10 pt,-\j*10 pt) {\pgfplotsretval};
    \fi
    };
};
\draw[->]        (-0.5,0)  -- node [above,midway,rotate=90] {Dimension} (-0.5,-3);
\draw[->]        (0,0.5)   -- node [above,midway] {Indices} (3,0.5);
\end{tikzpicture}
\caption{An abstract visualization of the first 40 indices $\boldk$ detected when using the Fourier series parametrization. The leftmost column contains the index $\boldk$ corresponding to the largest (in absolute value) basis coefficient $\hat{u}_{\boldk}$, the second column the index for the second largest and so on. The rows identify the $10$ dimensions corresponding to the variables $x$ and $a_{-4},\ldots ,a_{4}$ from top to bottom in this order. Zeros are neglected to preserve clarity.}
\label{fig:pois1a_ind}
\end{figure}
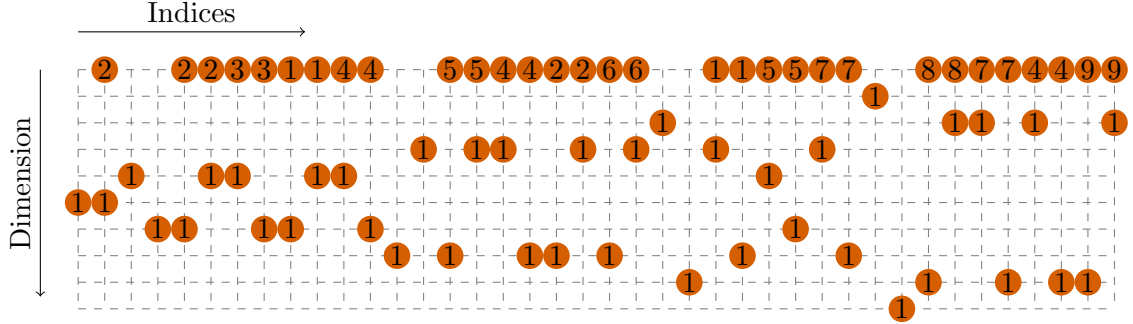

The detected indices show a clear structure as can be seen for the first indices in Figure~\ref{fig:pois1a_ind}. For all dimensions corresponding to an entry of the Fourier coefficient vector $\bolda$, there exists no other entry than $0$ or $1$. This effect is not caused by the particular sparsity $s$ we have chosen, since the algorithm already neglects every other possible entry (so the numbers from $2$ to $64$ for our choice of $N$) in the single component identification step, cf.\ Step 1 in Algorithm \ref{alg:main}, such that they can never appear at all in these dimensions. This result is exactly what we expected knowing the explicit formula \eqref{eq:u_pois_anal_trans}, since therein all the Fourier coefficients $a_{-4},\ldots,a_4$ appear only linearly. Additionally, the only indices where the entry corresponding to $a_0$ is non-zero are the first and second one, which can be also seen in Figure \ref{fig:pois1a_ind} as the first and second column. Again, this matches our expectations, since $a_0$ only appears in the first term of \eqref{eq:u_pois_anal_trans}, which can be rewritten in terms of the Chebyshev polynomials as
\begin{align*}
\frac{a_0}{8}(1-\tilde{x}^2) &= \frac1{8\sqrt2} T_1(a_0) \left(\frac12 T_0(x) - \frac1{2\sqrt2} T_2(x)\right) \prod_{\substack{\ell=-\frac{n-1}{2} \\ \ell\neq 0}}^{\frac{n-1}{2}} T_0(a_\ell) \\
&= \frac{1}{16\sqrt{2}} T_{\boldk^{(1)}}(x,\bolda) - \frac{1}{32} T_{\boldk^{(2)}}(x,\bolda),
\end{align*}
with
\begin{align*}
\boldk^{(1)} &= [\underbrace{0}_{T_0(x)},\underbrace{0,0,0,0}_{T_0(a_\ell)},\underbrace{1}_{T_1(a_0)},\underbrace{0,0,0,0}_{T_0(a_\ell)}]^T \\
\intertext{and} 
\boldk^{(2)} &= [\underbrace{2}_{T_2(x)},\underbrace{0,0,0,0}_{T_0(a_\ell)},\underbrace{1}_{T_1(a_0)},\underbrace{0,0,0,0}_{T_0(a_\ell)}]^T \\
\end{align*}
being the indices mentioned above. Finally, we would expect no couplings between the different Fourier coefficients $a_{-4},\ldots,a_4$ since they never appear together in the parts of the sum in the right part of \eqref{eq:u_pois_anal_trans}. While this behavior can be observed for the first detected indices in Figure \ref{fig:pois1a_ind}, this does not hold for all of our detected indices. At some point, the value of the remaining true Chebyshev coefficients become so small, that the algorithm can not distinguish their corresponding indices from false ones like $[2,1,0,1,1,1,0,1,0,1]^T$, which seem to produce similar coefficient values due to small numerical errors. However, since the size of the coefficients where this effect happens is already very small, i.e.\ about $10^{-8}$, this does not harm the overall approximation. Obviously, this minor problem is simply caused by the large sparsity $s=1000$ and could also be prevented up to some extend by using a search space $\Gamma$ that does not contain indices $\boldk$ with so many non-zero entries.

\begin{example}[High-dimensional extension of the detected index set $I$]
As stated already in Section \ref{sec:introduction}, we can use the structure of the detected index set $I$ to extend our approach to even higher dimensions. We demonstrate this approach here for the current differential equation with the Fourier series parametrization due to its simple structure and our explicit knowledge of the true solution \eqref{eq:u_pois_anal_trans}.

We increase the number $n$ of Fourier coefficients $a_\ell$ used to $99$ in order to achieve a better resolution of the right-hand side function $f$ than before. Obviously, this leads us to the approximation of the now $100$-dimensional function $u(x,\bolda)$. However, analyzing the detected index set $I$ from our $10$-dimensional test above, we can construct a good index set $I$ directly by generalizing the main structural features of $I$. In detail, we will construct our new index set $I$ in the following way: \begin{itemize}
\item The first dimension (corresponding to the spatial variable $x$) may contain any number from $0$ to $N_x$.
\item The entries of the dimensions $2$ to $100$ are either all zero or contain at most one non-zero entry. This non-zero entry, if existing, must be $1$.
\end{itemize}
This index set $I$ then contains $(N_x+1) \cdot 100$ indices of a similar structure as in Figure \ref{fig:pois1a_ind}. Note that we did not include the fact, that there were only two indices $\boldk^{(1)}$ and $\boldk^{(2)}$ with a non-zero entry in the dimension corresponding to $a_0$.

We perform our test using $N_x = 999$, so using an index set $I$ containing $10^5$ indices $\boldk$. We compute the corresponding basis cofficients $u_\boldk$ using the same Chebyshev multiple rank-1 lattice approach from \cite{Kae25} as before in Algorithm \ref{alg:main}. Our full dimension-incremental algorithm in just $10$ dimensions already needed around $400000$ samples for this simple example. Now we only need about $300000$ samples to approximate all the basis coefficients $\hat{u}_\boldk, \boldk \in I$ in this 100-dimensional example. This would be an impossible goal when applying our full algorithm \ref{alg:main} directly to the $100$-dimensional approximation problem instead. Especially for real applications, where the sampling values are not generated by an explicit formula but by a differential equation solver (like the FEM), the reduction of the amount of samples needed is an important tool, since the corresponding calls of the differential equation solver will be the dominating part of the computational complexity of the whole algorithm. The same problem appeared in \cite{KaPoTa22} and was the main motivation for the method proposed there. The relative approximation errors range from $10^{-8}$ to at most $10^{-6}$, which is a further improvement compared to the relative errors for the range $[-1,1]$ in Figure \ref{fig:pois1a_err}. Note that all these tests were performed with right-hand side functions $f$ of the form \eqref{eq:pois1d_f_four}, so as before there was no error in the discretization of this function. However, the higher resolution of this approach with $n=99$ allows for a much better discretization error (of the function $f$) if we work with more general right-hand side functions $f$.
\end{example}
 
\subsubsection{B-spline parametrization} \label{subsubsec:poisson_1d_bs}
We approximate the right-hand side $f$ by a sum of $n$ B-splines, i.e.
\begin{align*}
f(x) \approx \sum_{\ell=0}^{n-1} a_\ell B_\ell^{(m)}(x),
\end{align*}
where $B_\ell^{(m)}$ are versions of the cardinal B-spline $B^{(m)}$ of order $m$. Originally, they are recursively defined via
\begin{align*}
B^{(1)}(x) \coloneqq \begin{cases}
1, & -\frac12 < x < \frac12 \\
0, & \text{otherwise} 
\end{cases} \quad \text{and} \quad B^{(m)}(x) \coloneqq \int_{x-\frac12}^{x+\frac12} B^{(m-1)}(y)\diff y.
\end{align*}
We use the additional index $\ell$ to indicate, that we scaled and shifted them w.r.t.\ the interval $[0,1]$ and the desired amount of B-splines $n$, such that their peaks are equidistantly spaced along the interval and each spline overlaps $m-1$ neighboring splines in each direction.

This time, we use a classical differential equation solver to acquire the sample values $u(x^\ast,\bolda^\ast)$ for any sampling point $(x^\ast,\bolda^\ast)$. In particular, we will apply the function \texttt{solve\textunderscore bvp} from the submodule \texttt{scipy.integrate} here, which is capable of solving first order systems of ODEs with two-point boundary conditions. As mentioned in Remark \ref{rem:transform}, we need to transform the points $(x,\bolda)$ such that they fit to the domain $\D = [-1,1]^{n+1}$. We perform the same steps as in the previous example by transforming $\mathcal{T} x = \frac12(x+1)$ and restricting $\bolda \in [-1,1]^n$ throughout this example. Hence, once again the final function, which we are going to approximate, is $\tilde{u}(\tilde{x},\bolda) \coloneqq u(\mathcal{T}\tilde{x},\bolda) = u(\frac12(\tilde{x}+1),\bolda)$.

Additionally, we also use the same amount $n=9$ of Spline coefficients $a_\ell$, again resulting in $10$-dimensional approximation problem. The sparsity $s=1000$ and the extension $N=64$ of the search space $\Gamma$ also remain the same. With the differential equation solver needing roughly $0.1$ seconds per call, we ended up with an overall runtime of our method of almost $6$ hours.

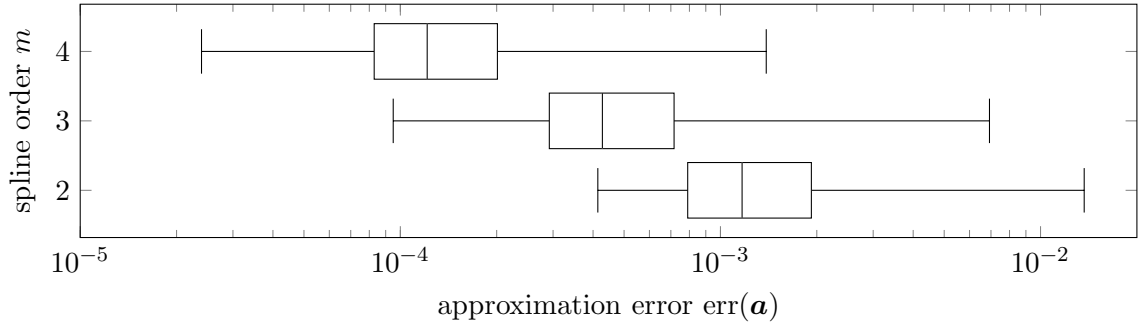
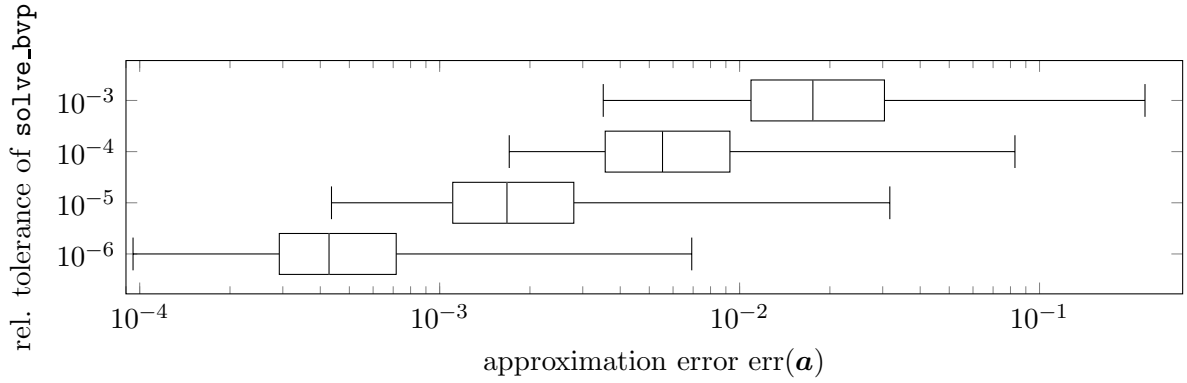
\begin{figure}[t]
\centering
\begin{subfigure}[c]{\textwidth}
\begin{tikzpicture}
\begin{axis}[
	boxplot/draw direction=x,
	xmode=log,
	xmin=0.00001,
	xmax=0.02,
	ytick={1,2,3},
	yticklabels={$2$,$3$,$4$},
	xlabel={approximation error $\text{err}(\bolda)$},
	ylabel={spline order $m$},
	width=\textwidth,
    height=0.30\textwidth,
	]
	\pgfplotstableread[col sep=comma]{1b_err.csv}\datatable;
	\addplot+[
  		boxplot prepared from table={
    			table=\datatable,
    			row=0,
    			lower whisker=lw,
    			upper whisker=uw,
    			lower quartile=lq,
    			upper quartile=uq,
    			median=med
  		}, my boxplot style,
    ] coordinates {};
    \addplot+[
  		boxplot prepared from table={
    			table=\datatable,
    			row=1,
    			lower whisker=lw,
    			upper whisker=uw,
    			lower quartile=lq,
    			upper quartile=uq,
    			median=med
  		}, my boxplot style,
    ] coordinates {};
    \addplot+[
  		boxplot prepared from table={
    			table=\datatable,
    			row=2,
    			lower whisker=lw,
    			upper whisker=uw,
    			lower quartile=lq,
    			upper quartile=uq,
    			median=med
  		}, my boxplot style,
    ] coordinates {};
    \end{axis}
\end{tikzpicture}
\caption{Varying the spline order $m$ for fixed relative tolerance $10^{-6}$ (default).}
\label{fig:pois1b_err_m}
\end{subfigure}
~
\begin{subfigure}[c]{\textwidth}
\begin{tikzpicture}
\begin{axis}[
	boxplot/draw direction=x,
	xmode=log,
	xmin=0.00009,
	xmax=0.3,
	ytick={4,3,2,1},
	yticklabels={$10^{-3}$,$10^{-4}$,$10^{-5}$, $10^{-6}$},
	xlabel={approximation error $\text{err}(\bolda)$},
	ylabel={rel. tolerance of \texttt{solve\textunderscore bvp}},
	width=\textwidth,
    height=0.30\textwidth,
	]
	\pgfplotstableread[col sep=comma]{1b_err.csv}\datatable;
	\addplot+[
  		boxplot prepared from table={
    			table=\datatable,
    			row=1,
    			lower whisker=lw,
    			upper whisker=uw,
    			lower quartile=lq,
    			upper quartile=uq,
    			median=med
  		}, my boxplot style,
    ] coordinates {};
	\addplot+[
  		boxplot prepared from table={
    			table=\datatable,
    			row=5,
    			lower whisker=lw,
    			upper whisker=uw,
    			lower quartile=lq,
    			upper quartile=uq,
    			median=med
  		}, my boxplot style,
    ] coordinates {};
    \addplot+[
  		boxplot prepared from table={
    			table=\datatable,
    			row=4,
    			lower whisker=lw,
    			upper whisker=uw,
    			lower quartile=lq,
    			upper quartile=uq,
    			median=med
  		}, my boxplot style,
    ] coordinates {};
    \addplot+[
  		boxplot prepared from table={
    			table=\datatable,
    			row=3,
    			lower whisker=lw,
    			upper whisker=uw,
    			lower quartile=lq,
    			upper quartile=uq,
    			median=med
  		}, my boxplot style,
    ] coordinates {};
    \end{axis}
\end{tikzpicture}
\caption{Varying the relative tolerance for fixed spline order $m=3$.}
\label{fig:pois1b_err_acc}
\end{subfigure}
\caption{The relative approximation error $\text{err}(\bolda)$ for $10000$ randomly drawn $\bolda$ for different choices of the spline order $m$ and the relative tolerance of the solver function \texttt{solve\textunderscore bvp}. The box-and-whisker plots show the median, the first and the second quartile as well as the maximal and minimal error observed.}
\label{fig:pois1b_err}
\end{figure}

The approximation error shown in Figure \ref{fig:pois1b_err} is derived by evaluating our approximation as well as the solution given by the function \texttt{solve\textunderscore bvp} on $1000$ equidistant spatial points and considering the respective $\ell_2$-error $\text{err}(\bolda)$. The reference solution for the error erstimation was obtained by using a relative tolerance of $10^{-9}$ in \texttt{solve\textunderscore bvp}. As in the previous example, we use box-and-whisker plots to visualize the statistical distribution of the results, i.e., the range and the median of the observed errors as well as their quartiles.

Figure \ref{fig:pois1b_err_m} shows the results for different choices of the spline order $m$ when parametrizing the right-hand side function $f$. The piecewise linear splines ($m=2$) result in rather unsatisfying errors, probably caused by the lack of smoothness. Higher-order splines as $m=3$ and $m=4$ provide better results, especially when investigating the range and the worst case of the possible errors $\text{err}(\bolda)$. Although the overall error size in Figure \ref{fig:pois1b_err_m} might seem a bit large, it is matching the default relative tolerance $10^{-6}$ of the function \texttt{solve\textunderscore bvp}, which we used for these tests. A lower accuracy of the underlying differential equation solver leads to a lower accuracy of our method, which can be observed in Figure \ref{fig:pois1b_err_acc}. Here, we fixed the spline order $m=3$ and varied the relative tolerance of the function \texttt{solve\textunderscore bvp}. As mentioned before, we will not go into further detail about the properties of the differential equation solvers used. However, we wanted to briefly mention the influence of the accuracy of the underlying solver at least for this first example.

\begin{figure}[t]
\centering
\begin{subfigure}[c]{\textwidth}
\begin{tikzpicture}
\pgfplotstableread{ind.dat}\mytable
\foreach \i in {0,...,39} {
    \draw[dashed, thin, gray, opacity=0.5] (\i*10 pt,0) -- (\i*10 pt,-90 pt);
}
\foreach \j in {0,...,9} {
    \draw[dashed, thin, gray, opacity=0.5] (0,-\j*10 pt) -- (390 pt,-\j*10 pt);
}
\foreach \i in {40,...,79}{
    \foreach \j in {0,...,9}{
    \pgfplotstablegetelem{\i}{\j}\of\mytable
    \ifnum\pgfplotsretval=0\relax\else
    \node[circle, minimum size=10pt, inner sep=0pt, fill=color6, opacity=0.\pgfplotsretval, text opacity= 0.75] at (\i*10-400 pt,-\j*10 pt) {\pgfplotsretval};
    \fi
    };
};
\draw[->]        (-0.5,0)  -- node [above,midway,rotate=90] {Dimension} (-0.5,-3);
\draw[->]        (0,0.5)   -- node [above,midway] {Indices} (3,0.5);
\end{tikzpicture}
\caption{The detected indices from number $1$ to number $40$.}
\label{fig:pois1b_ind}
\end{subfigure}
~
\begin{subfigure}[c]{\textwidth}
\centering
\begin{tikzpicture}
\pgfplotstableread{ind.dat}\mytable
\foreach \i in {0,...,39} {
    \draw[dashed, thin, gray, opacity=0.5] (\i*10 pt,0) -- (\i*10 pt,-90 pt);
}
\foreach \j in {0,...,9} {
    \draw[dashed, thin, gray, opacity=0.5] (0,-\j*10 pt) -- (390 pt,-\j*10 pt);
}
\foreach \i in {80,...,119}{
    \foreach \j in {0,...,9}{
    \pgfplotstablegetelem{\i}{\j}\of\mytable
    \ifnum\pgfplotsretval=0\relax\else
    \node[circle, minimum size=10pt, inner sep=0pt, fill=color6, opacity=0.1*\pgfplotsretval, text opacity= 0.75] at (\i*10-800 pt,-\j*10 pt) {\pgfplotsretval};
    \fi
    };
};
\draw[->]        (-0.5,0)  -- node [above,midway,rotate=90] {Dimension} (-0.5,-3);
\draw[->]        (0,0.5)   -- node [above,midway] {Indices} (3,0.5);
\end{tikzpicture}
\caption{The detected indices from number $201$ to number $240$.}
\label{fig:pois1b2_ind}
\end{subfigure}
\caption{Abstract visualizations of $40$ detected indices $\boldk$ (from left to right) for Example \ref{subsubsec:poisson_1d_bs}. The indices $\boldk$ are sorted in descending order according to the size of the corresponding approximated coefficient $\hat{u}_{\boldk}$. The rows identify the $10$ dimensions corresponding to the variables $x$ and $a_{-4},\ldots ,a_{4}$ from top to bottom in this order. Zeros are neglected to preserve clarity.}
\end{figure}

Figure \ref{fig:pois1b_ind} shows two parts of the detected index set $I$ for the spline order $m=3$. As in Example \ref{subsubsec:poisson_1d_fc}, we notice a sparse structure of the first detected indices. This time, we do not have an explicit representation of the true solution $u$ and only use approximations of the solution given by the differential equation solver as our samples. Hence, we are not capable of comparing these indices and the corresponding values to the true ones as before. However, the structure, which can be observed in Figure \ref{fig:pois1b_ind}, is still highly reasonable: It shows, that the algorithm is prioritizing two-dimensional couplings with small entries in the first dimension corresponding to the spatial variable $x$. For later indices as shown for example in Figure \ref{fig:pois1b2_ind}, there appear some higher-dimensional couplings and even some values greater than $1$ outside of the spatial dimension. The coupling B-spline coefficients $a_\ell$ are always adjacent. This is caused by the overlapping nature of the B-splines. Even for later indices (apart from numerical errors as described below) this behavior will continue.

On the other hand, each of the $1000$ detected indices contains at least one non-zero entry in the dimensions $2$ to $10$, i.e., the corresponding Chebyshev series does not contain a single term that depends only on $x$. Unfortunately, there are also some artifact indices again, which do not contain a single zero entry but unreasonably large numbers ($\geq 10$) in these dimensions. As before, these are mainly caused by numerical errors and could be reduced by choosing the search space $\Gamma$ more restrictively.

Finally, all the computed coefficients are real-valued this time. While the domain of the coefficients $\bolda$ is the same as before, multiplying them with the cardinal B-splines instead of Fourier terms causes the source function $f$ and therefore also the solution $u$ to be real-valued for each possible coefficient $\bolda$. 

\subsection{A piece-wise continuous differential equation}
\label{subsec:pw_ode}
As a second one-dimensional example, we consider the ordinary differential equation
\begin{equation} \label{eq:pwc_ode}
\begin{aligned} 
-\frac{d}{d x}(a(x) \frac{d}{d x} u(x)) &= f(x), & x \in (-1,1), \\
u(-1) = u(1) &= 0, &
\end{aligned}
\end{equation}
with the piece-wise constant coefficient function
\begin{align*}
a(x) = \begin{cases}
\frac12, & x\in (-1,0),\\
1, & x\in [0,1).
\end{cases}
\end{align*}
This example was investigated in \cite[Sec.~2.3]{LuZh23} and threw up tremendous problems when using physics-informed neural networks (PINNs) since it has no classical but only a weak solution $u$ for the given right-hand side function $f$
\begin{align}\label{eq:pwc_f}
f(x) = \begin{cases}
0, & x\in (-1,0),\\
-2, & x\in [0,1).
\end{cases}
\end{align}
Therefore, we are interested in solving \eqref{eq:pwc_ode} using our approach and comparing the result afterwards for this particular right-hand side function $f$. The exact solution for this scenario is also given in \cite{LuZh23} and reads as
\begin{align}\label{eq:pwc_ex}
u(x) = \begin{cases}
-\frac23 x - \frac23, & x\in (-1,0),\\
x^2 - \frac13 x - \frac23, & x\in [0,1).
\end{cases}
\end{align}

For this ODE, we have the differential operator $\L = -\frac{d}{d x} a(x) \frac{d}{d x}$ and choose $\U = H_0^1((-1,1))$ and $\F = H^{-1}((-1,1))$ as the function spaces as well as the approximation domain $\D = [-1,1]^{n+1}$ and the tensorized Chebyshev polynomials as the BOPB. In order to resolve \eqref{eq:pwc_f} properly, we choose a discretization of $f$ similar to Section \ref{subsubsec:poisson_1d_bs} using B-splines of order $m=1$, so characteristic functions on non-overlapping intervals. Precisely, we resolve the right-hand side $f$ as
\begin{align}\label{eq:pwc_char}
f(x) = \sum_{\ell=0}^7 b_\ell \mathds{1}_{[-1+\frac{\ell}{4},-1+\frac{\ell+1}{4}]}(x),
\end{align}
such that the particular function $f$ given in \eqref{eq:pwc_f} is obtained exactly for the spline coefficients $\boldb = [0,0,0,0,-2,-2,-2,-2]^T$. Then, the general exact solution reads as 
\begin{align}\label{eq:pwc_gen}
u(x,\boldb) = \begin{cases}
\sum_{\ell=0}^7 -2b_\ell W_\ell(x) + 2 C_1 x + C_2, & x \in (-1,0), \\
\sum_{\ell=0}^7 -b_\ell W_\ell(x) + C_1 x + C_2, & x \in [0,1), \\
\end{cases}
\end{align}
with
\begin{align*}
W_\ell(x) \coloneqq \begin{cases}
0, & x \in (-1, -1+\frac{\ell}4), \\
\frac12 x^2 + (1-\frac{\ell}4)x + \frac12(-1+\frac{\ell}4)^2, & x \in [-1+\frac{\ell}4,-1+\frac{\ell+1}4), \\
\frac14 x + \frac7{32}-\frac{\ell}{16}, & x \in [-1+\frac{\ell+1}4,1)\\
\end{cases}
\end{align*}
being the second anti-derivative of the characteristic function $\mathds{1}_{[-1+\frac{\ell}{4},-1+\frac{\ell+1}{4}]}(x)$. The boundary conditions from \eqref{eq:pwc_ode} yield $C_2=2C_1$ and $C_1 = \sum_{\ell=0}^7 b_\ell \frac{15-2\ell}{96}$.

To obtain the necessary samples of the (weak) solution of \eqref{eq:pwc_ode}, we utilize the popular open-source computing platform FEniCS and its Python interface. The underlying finite element method is very well capable of computing the weak solution of \eqref{eq:pwc_ode}, making it a perfect tool for our black-box sampling step. We use a rather coarse mesh with only $100$ nodes. Note that the domain of the spatial variable $x$ is already $[-1,1]$ this time and needs no further transformation. On the other hand, we scale the spline coefficients $b_\ell$ with a factor of $2$ (or $\frac12$, respectively) to cover the range $[-2,2]^8$, such that the particular right-hand side $f$ given in \eqref{eq:pwc_f} can be resolved exactly as described above.

The amount $n=8$ of spline coefficients $b_\ell$ is already fixed such that we have the overall dimension $d=9$ for this problem. We use the sparsity $s=4000$ with the extension $N=256$ for this test example. The dimension-incremental method needed about $75$ minutes for this approach. FEniCS itself computes one solution sample for this problem in about $10^{-2}$ seconds.

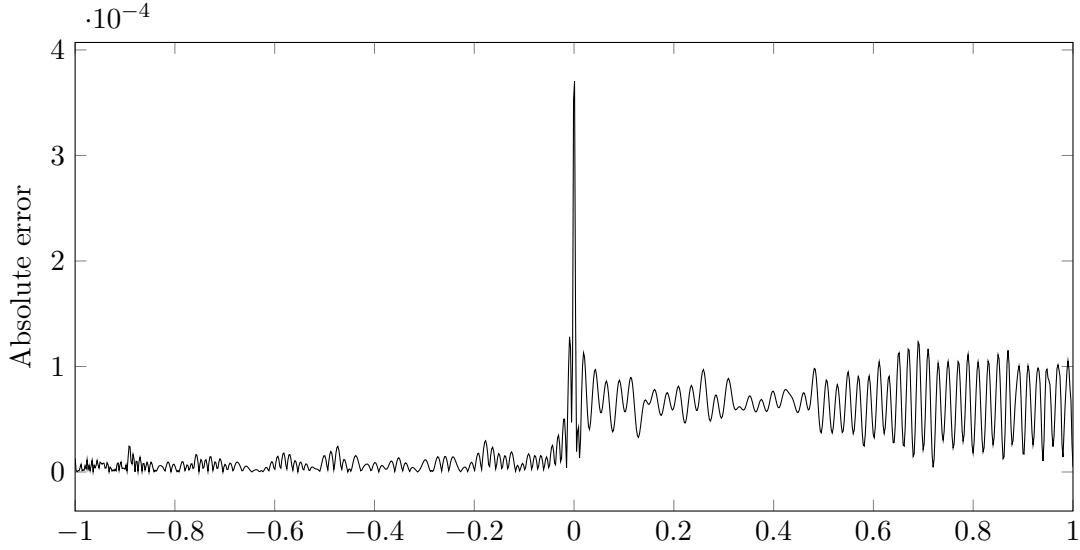
\begin{figure}[t]
\centering
\begin{tikzpicture}
 \begin{axis}[
 xmin = -1,
 xmax = 1,
 ylabel={Absolute error},
 width=0.95\textwidth,
 height=0.50\textwidth,]
 \addplot[solid, color=black] table[x index=0, y index=1] {2_error_data.txt};
 \end{axis}
\end{tikzpicture}
\caption{The (absolute) pointwise approximation error of our approximation when using the right-hand side \eqref{eq:pwc_f} compared to the exact solution \eqref{eq:pwc_ex}.}
\label{fig:pwc_err}
\end{figure}

Figure \ref{fig:pwc_err} illustrates the pointwise error of our approximation for the particular function $f$ given in \eqref{eq:pwc_f} when compared to the true solution \eqref{eq:pwc_ex}. We note that the kink of the exact solution $u$ at the point $x=0$ leads to larger pointwise errors in this region, which is not surprising given our smooth basis functions and hence the smoothness of our approximation. Further, it appears that in the right half of the domain, the pointwise errors are no longer unbiased but oscilliate around some constant greater than zero.

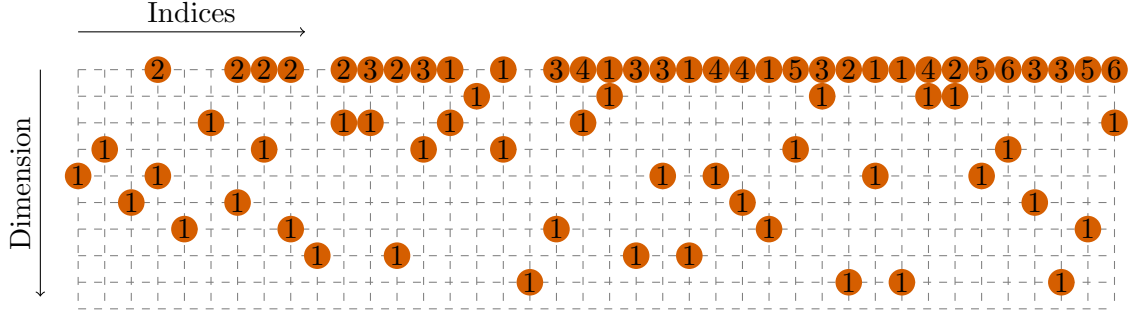
\begin{figure}[t]
\centering
\begin{tikzpicture}
\pgfplotstableread{ind.dat}\mytable
\foreach \i in {0,...,39} {
    \draw[dashed, thin, gray, opacity=0.5] (\i*10 pt,0) -- (\i*10 pt,-90 pt);
}
\foreach \j in {0,...,9} {
    \draw[dashed, thin, gray, opacity=0.5] (0,-\j*10 pt) -- (390 pt,-\j*10 pt);
}
\foreach \i in {120,...,159}{
    \foreach \j in {0,...,8}{
    \pgfplotstablegetelem{\i}{\j}\of\mytable
    \ifnum\pgfplotsretval=0\relax\else
    \node[circle, minimum size=10pt, inner sep=0pt, fill=color6, opacity=0.1*\pgfplotsretval, text opacity= 0.75] at (\i*10-1200 pt,-\j*10 pt) {\pgfplotsretval};
    \fi
    };
};
\draw[->]        (-0.5,0)  -- node [above,midway,rotate=90] {Dimension} (-0.5,-3);
\draw[->]        (0,0.5)   -- node [above,midway] {Indices} (3,0.5);
\end{tikzpicture}
\caption{An abstract visualization of the first 40 indices $\boldk$ detected for the piece-wise continuous differential equation example. The indices $\boldk$ are ordered by the absolute values of their corresponding approximated coefficients $\hat{u}_\boldk$ in descending order from left to right. The rows identify the $9$ dimensions corresponding to the spatial variable $x$ and the $n=8$ spline coefficients $\boldb$ used. Zeros are neglected to preserve clarity.}
\label{fig:pwc_ind}
\end{figure}

The detected index set, partially illustrated in Figure \ref{fig:pwc_ind}, shows the usual structure from the previous example. The linear dependence of the spline coefficients $\boldb$ in the analytical solution \eqref{eq:pwc_gen} is detected perfectly, neglecting any entries larger than $1$ in these dimensions already in Step 1 of Algorithm \ref{alg:main}. This is highly remarkable when compared to our first example in Section \ref{subsubsec:poisson_1d_fc}, since we are using a numerical solver instead of an analytical solution this time. So while the samples of the solution $u$ contain small numerical errors, our algorithm still successfully neglects the unnecessary values in Step 1. The first dimension, corresponding to the spatial variable $x$, contains again larger values, caused by the highly piecewise structure of the solution \eqref{eq:pwc_gen}.

\subsection{The multi-dimensional Poisson equation}
\label{subsec:poisson_hd}

While the previous examples considered ordinary differential equations, we now progress to partial differential equations with the two-dimensional version of \eqref{eq:poisson_1d}. The general Poisson equation with homogeneous Dirichlet boundary conditions is given by
\begin{equation} \label{eq:poisson_hd}
\begin{aligned} 
-\Delta u(\boldx) &= f(\boldx), & \boldx \in \Omega, \\
u(\boldx) &= 0, & \boldx \in \delta\Omega,
\end{aligned}
\end{equation}
with the spatial domain $\Omega = (0,1)^d$. We restrict ourself to the two-dimensional version $d=2$ in this work, while $d=3$ is also a common setting in applications. With the differential operator $\L = -\Delta$ and the homogeneous Dirichlet boundary conditions, we use the common function spaces $\U = H_0^1(\Omega)$ and $\F = H^{-1}(\Omega)$, which is the direct generalization of the function spaces used in Section \ref{subsec:poisson_1d}. Also, the approximation space $\H = L_2(\D)$ is again equipped with the tensorized Chebyshev polynomials on the now $n+2$-dimensional domain $\D = [-1,1]^{n+2}$.

Motivated by the example from Section \ref{subsec:poisson_1d}, we use a two-dimensional Fourier series to parametrize the right-hand side function $f$ in this example. In detail, we parametrize $f$ by
\begin{align*}
f(\boldx) \approx \sum_{\boldell \in J} a_\boldell \e^{2\pi\ii \boldell \boldx}
\end{align*}
with the index set $J$, again containing a total of $n$ indices. As in the one-dimensional case, this choice should result in a rather simple detected index set $I$ when applying our algorithm.

As in Section \ref{subsec:pw_ode}, we use FEniCS to solve the PDE \eqref{eq:poisson_hd} with the finite element method for given $\bolda$ this time. For the finite element mesh, we used a uniform unit square mesh with $51$ equidistant points in each direction. Since each square cell is split into two triangular elements, we thus end up with $5000$ finite elements for our approximation. As in the ODE examples, we still need to transform the sampling points $(\boldx,\bolda)$. Hence, we proceed similarly as in Section \ref{subsubsec:poisson_1d_fc} by using the transformation $\mathcal{T}\boldx = \frac12(\boldx + \boldone)$ and simply restricting $\bolda \in [-1,1]^n$.

We set $J \coloneqq \{-1,0,1\}^2$ in order to have $n=9$ Fourier coefficients. Combined with the spatial dimension $d=2$, we end up with an 11-dimensional approximation problem this time. Further, we choose the sparsity $s=1000$ and the extension $N=64$ of the search space $\Gamma$. FEniCS needs up to $3$ seconds for each single solution this time, leaving our dimension-incremental method with a total runtime of roughly $100$ minutes.

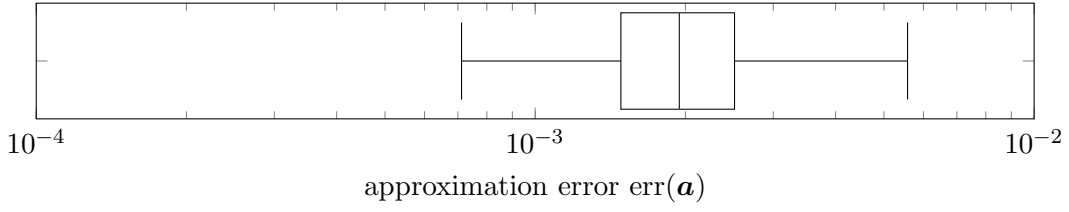
\begin{figure}[t]
\centering
\begin{tikzpicture}
\begin{axis}[
	boxplot/draw direction=x,
	xmode=log,
	xmin=0.0001,
	xmax=0.01,
	ytick={1},
	yticklabels={},
	xlabel={approximation error $\text{err}(\bolda)$},
	width=0.95\textwidth,
    height=0.20\textwidth,
	]
	\pgfplotstableread[col sep=comma]{3_err.csv}\datatable;
	\addplot+[
  		boxplot prepared from table={
    			table=\datatable,
    			row=0,
    			lower whisker=lw,
    			upper whisker=uw,
    			lower quartile=lq,
    			upper quartile=uq,
    			median=med
  		}, my boxplot style,
    ] coordinates {};
    \end{axis}
\end{tikzpicture}
\caption{The relative approximation error $\text{err}(\bolda)$ for $10000$ randomly drawn $\bolda$ for the two-dimensional Poisson equation example. The box-and-whisker plots show the median, the first and the second quartile as well as the maximal and minimal error observed.}
\label{fig:pois2_err}
\end{figure}

Figure \ref{fig:pois2_err} illustrates the relative approximation error $\text{err}(\bolda)$ as before in the one-dimensional case using $10000$ randomly drawn coefficients $\bolda$. Note that this error is computed by comparing our approximation to the solution the FEM solver produces for the given $\bolda$ on the nodes of the FE mesh using the same parameters as we did during the execution of our algorithm. We observe errors of sizes around $10^{-3}$, which are obviously larger than before in Section \ref{subsubsec:poisson_1d_fc}. However, this effect is primarily caused by the fact, that we are no longer using a direct representation of the analytic solution as for the one-dimensional example to generate our sampling points. The error sizes are still reasonable and can compete with the used PDE solver, taking into account the sparsity $s$ and extension $N$ used here.

\begin{figure}[t]
\centering
\begin{tikzpicture}
\pgfplotstableread{ind.dat}\mytable
\foreach \i in {0,...,39} {
    \draw[dashed, thin, gray, opacity=0.5] (\i*10 pt,0) -- (\i*10 pt,-100 pt);
}
\foreach \j in {0,...,10} {
    \draw[dashed, thin, gray, opacity=0.5] (0,-\j*10 pt) -- (390 pt,-\j*10 pt);
}
\foreach \i in {160,...,199}{
    \foreach \j in {0,...,10}{
    \pgfplotstablegetelem{\i}{\j}\of\mytable
    \ifnum\pgfplotsretval=0\relax\else
    \node[circle, minimum size=10pt, inner sep=0pt, fill=color6, opacity=0.\pgfplotsretval, text opacity= 0.75] at (\i*10-1600 pt,-\j*10 pt) {\pgfplotsretval};
    \fi
    };
};
\draw[->]        (-0.5,0)  -- node [above,midway,rotate=90] {Dimension} (-0.5,-3);
\draw[->]        (0,0.5)   -- node [above,midway] {Indices} (3,0.5);
\end{tikzpicture}
\caption{An abstract visualization of the first 40 indices $\boldk$ detected for the two-dimensional Poisson equation example. The indices $\boldk$ are ordered by the absolute values of their corresponding approximated coefficients $\hat{u}_\boldk$ in descending order from left to right. The rows identify the $11$ dimensions corresponding to the two spatial variables $\boldx=(x_1,x_2)^T$ and the $n=9$ Fourier coefficients $\bolda$ used. Zeros are neglected to preserve clarity.}
\label{fig:pois2_ind}
\end{figure}
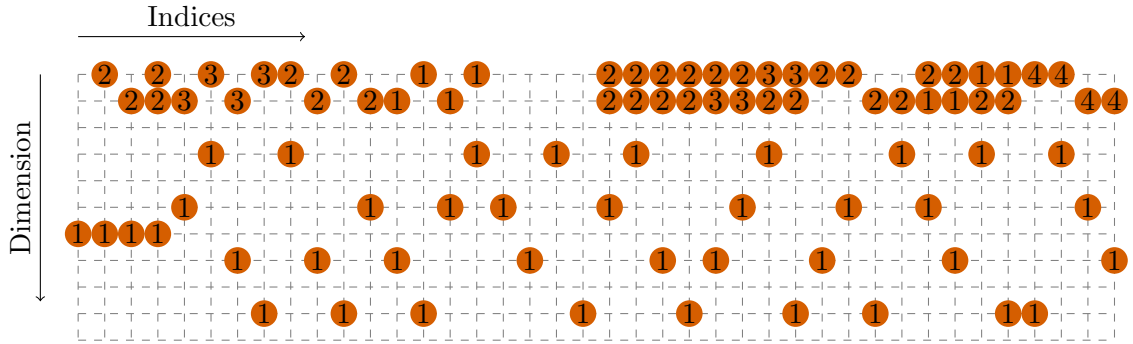

The structure of the detected index set $I$, where the first part is shown in Figure \ref{fig:pois2_ind}, is pretty similar to the one seen in Section \ref{subsubsec:poisson_1d_fc}. Even though we are not using the exact solution for our training samples anymore, our algorithm is still able to identify that for the dimensions corresponding to the Fourier coefficients $\bolda$ the only necessary entries are $0$ and $1$. The entries of the first two dimensions, corresponding to the spatial dimensions $\boldx$, also contain larger numbers, but are growing significantly slower than in the one-dimensional example. This is due to the fact, that in this two-dimensional case also all possible combinations of the entries in these two dimensions have to be exploited. Overall, the discovered structure resembles the one from our first example quite nicely and seems like the kind of structure for such an index set we would expect as the canonical generalization to multi-dimensional examples even without examining an analytical solution.

\subsection{A diffusion equation with an affine random coefficient}
\label{subsec:rand_diff_eq}

The differential equation \eqref{eq:pwc_ode} is also a one-dimensional diffusion equation, for which the coefficient $a$ could be also called diffusion coefficient. Now, we investigate a two-dimensional diffusion equation on $\Omega=[0,1]^2$ with a randomized diffusion coefficient $a$, which is not only a PDE instead of an ODE but also a parametrized differential equation, i.e.
\begin{equation} \label{eq:rand_diff_eq}
\begin{aligned} 
-\nabla\cdot(a(\boldx,\boldy) \nabla u(\boldx,\boldy)) &= f(\boldx), & \boldx \in \Omega,\, \boldy \in \Omega_\boldy, \\
u(\boldx,\boldy) &= 0, & \boldx \in \partial\Omega,\, \boldy \in \Omega_\boldy.
\end{aligned}
\end{equation}
Here, the differential operator $\nabla$ is always used w.r.t.\ the spatial variable $\boldx$. While there exist multiple kinds of randomized diffusion coefficients $a$, as can be seen for example in \cite{KaPoTa22}, we will only work with an affine random coefficient $a$ here. In more detail, we consider the particular example from \cite[Sec.~11]{EiGiSchwZa14}, where we have for $n_\boldy=20$ the affine coefficient
\begin{align*}
a(\boldx,\boldy) \coloneqq 1 + \sum_{j=1}^{n_\boldy} y_j \psi_j(\boldx), & & \boldx \in \Omega,\, \boldy \in [-1,1]^{20}
\end{align*}
with the random variables $\boldy \sim \mathcal{U}([-1,1]^{n_\boldy})$ and
\begin{align*}
&\psi_j(\boldx) \coloneqq c j^{-\mu} \cos(2\pi m_1(j)\,x_1)\, \cos(2\pi m_2(j)\,x_2), & \boldx \in \Omega,\, j\geq 1.
\end{align*}
Here, $c > 0$ is a constant and $\mu > 1$ the decay rate. In our numerical example below, we use the values $c=0.9/\zeta(2)$ and $\mu = 2$ already used in \cite{EiGiSchwZa14}. Further, $m_1(j)$ and $m_2(j)$ are defined as
\begin{align*}
m_1(j) \coloneqq j-\frac{k(j) (k(j)+1)}{2} \quad \text{and} \quad m_2(j) \coloneqq k(j)-m_1(j)\\
\end{align*}
with $k(j) \coloneqq \lfloor -1/2 + \sqrt{1/4 + 2j} \rfloor$. For some explicit values of $m_1(j), m_2(j)$ and $k(j)$ as well as more details on this differential problem, see \cite{EiGiSchwZa14}. As before, we consider the common function spaces $\U = H_0^1(\Omega)$ and $\F = H^{-1}(\Omega)$ for this differential operator.

\begin{remark}
We already considered the numerical solution of this problem in \cite[Sec.~4.3]{KaPoTa22} using a slightly different approach. Therein, we discretize the spatial domain $\Omega = [0,1]^2$ and compute approximations like \eqref{eq:param_u_approx} with $\boldtheta\coloneqq\boldy$ in the Fourier setting for every fixed node $\boldx_g, g=1,\ldots,G$. The key ingredient there is, that the a priori unknown index set $I$ is chosen similarly for each of the $G$ approximations $S_I^\A u(\boldx_g,\cdot)$, which allows us to compute all these approximations using only a single call of a modification of the sparse approximation algorithm with slightly more samples and computation time needed. For a given random coefficient $\boldy^\ast$, we then compute the values of $S_I^\A u(\boldx_g,\boldy^\ast)$ at all the nodes $\boldx_g$ and interpolated between them to receive a solution on the complete domain $\Omega$.
\end{remark}

In contrast to all other examples considered in this work, we decided to use the fixed right-hand side $f \equiv 1$ without parametrization, since we are mainly interested in a comparison with the results from \cite{KaPoTa22}. Hence, we neglect the space $\F$ for this example and proceed with the solution operator $\G: \Omega_\boldy \rightarrow \U$ this time. The approximation space is still $\H = L_2(\D)$, again using the tensorized Chebyshev polynomials and $\D = [-1,1]^{n_\boldy + 2}$. The random variables $\boldy$ already match that domain, so we only transform $\boldx$ as usual using the transformation $\mathcal{T}\boldx = \frac12(\boldx + \boldone)$. Note that one could still use a similar approach as in the previous examples to easily parametrize the right-hand side $f$ as well in this example. As in Section \ref{subsec:poisson_hd}, we utilize FEniCS to solve the differential equation using $5000$ finite elements.

Since there is no parametrization of the right-hand side $f$, but $2$ spatial dimension as well as $n_\boldy = 20$ dimensions for the random variable $\boldy$, we still end up with a total of $22$ dimensions for our approximation problem. Unfortunately, the sampling complexity as well as the computational complexity of our approach using the Chebyshev basis include an exponential factor on the maximal number of non-zero entries of the indices $\boldk$ appearing in the candidate sets $K$ during step 2 of Algorithm \ref{alg:main}. In order to prevent cases, where numerical errors cause candidates $\boldk$ with (almost) non zeros in any dimensions, we impose a superposition dimension $d_s=7$ on our $22$-dimensional search space $\Gamma$ with extension $N=64$ for this example. Unfortunately, the more difficult structure of the differential equation lead up to a FEniCS computation time of up to $14$ seconds per call, resulting in a total runtime of about $60$ hours (for sparsity $s=1000$).

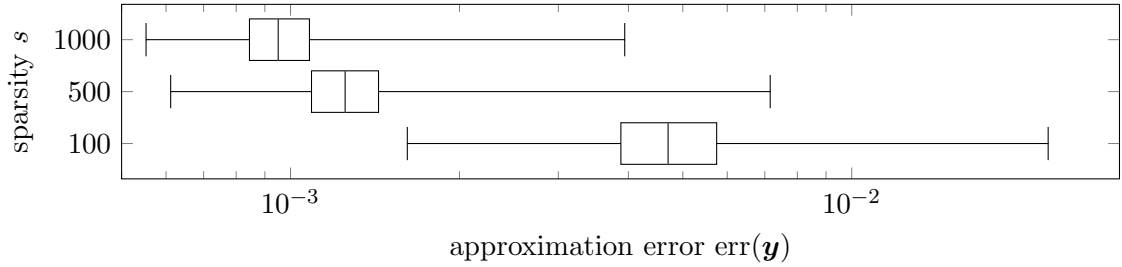
\begin{figure}[t]
\centering
\begin{tikzpicture}
\begin{axis}[
	boxplot/draw direction=x,
	xmode=log,
	xmin=0.0005,
	xmax=0.03,
	ytick={1,2,3},
	yticklabels={$100$, $500$, $1000$},
	xlabel={approximation error $\text{err}(\boldy)$},
	ylabel={sparsity $s$},
	width=0.95\textwidth,
    height=0.25\textwidth,
	]
	\pgfplotstableread[col sep=comma]{4_err.csv}\datatable;
	\addplot+[
  		boxplot prepared from table={
    			table=\datatable,
    			row=0,
    			lower whisker=lw,
    			upper whisker=uw,
    			lower quartile=lq,
    			upper quartile=uq,
    			median=med
  		}, my boxplot style,
    ] coordinates {};
    \addplot+[
  		boxplot prepared from table={
    			table=\datatable,
    			row=1,
    			lower whisker=lw,
    			upper whisker=uw,
    			lower quartile=lq,
    			upper quartile=uq,
    			median=med
  		}, my boxplot style,
    ] coordinates {};
    \addplot+[
  		boxplot prepared from table={
    			table=\datatable,
    			row=2,
    			lower whisker=lw,
    			upper whisker=uw,
    			lower quartile=lq,
    			upper quartile=uq,
    			median=med
  		}, my boxplot style,
    ] coordinates {};
    \end{axis}
\end{tikzpicture}
\caption{The relative approximation error $\text{err}(\boldy)$ for $10000$ randomly drawn variables $\boldy$ for the diffusion equation example for different sparsities $s$. The box-and-whisker plots show the median, the first and the second quartile as well as the maximal and minimal error observed.}
\label{fig:diff_pde_err}
\end{figure}

Figure \ref{fig:diff_pde_err} illustrates the relative approximation error $\text{err}(\boldy)$ for different sparsities $s$, this time with respect to the random variable $\boldy$. As in the previous example, the error is computed using the solution of the FEM with the same settings as comparison values. The error is again of reasonable size, even though this differential problem is significantly more difficult than the previous two-dimensional example in Section \ref{subsec:poisson_hd}. Further, the largest nodal error as considered in \cite[Sec.~4, Fig.~6]{KaPoTa22}, so the largest error at any of the nodes of the FE mesh when evaluating the approximation for $10000$ randomly drawn $\boldy$ and considering the respective $\ell_2$ norm, is just slightly larger than for the uniform sparse FFT from \cite{KaPoTa22}. This small increase is probably caused by the fact, that we are no longer focusing on particular nodes and basis expansions of the solution $u$ at these nodes, but a full basis expansion of $u$ also considering the spatial variable.

The structure of the $22$-dimensional index set is pretty similar to the previous examples and not illustrated here due to the high amount of dimensions. Surprisingly, the range of the entries in the dimensions corresponding to the random variables $\boldy$ is rather restricted. Already in the one-dimensional detections (Step 1 in Algorithm \ref{alg:main}), the algorithm does not detect a full range of $65$ possible entries (from $0$ to $N=64$), but less than $20$ possible entries for $y_1$, decaying rapidly down to only $4$ possible entries (so $0$, $1$, $2$ and $3$) for the later dimensions like $y_{15}$. While we already saw the extreme version of this behavior for the Poisson equation using the Fourier series parametrization, where the only possible entries were $0$ and $1$, we did not observe anything similar for the other examples like in Section \ref{subsubsec:poisson_1d_bs}.

\subsection{Heat equation}
\label{subsec:heat_eq}

Our next example is the heat equation in one dimension with homogeneous boundary conditions on the domain $\Omega = (0,L)$, i.e.
\begin{equation} \label{eq:heat_1d}
\begin{aligned} 
\partial_\tau u - \alpha^2 \partial_{xx} u &= 0, & x \in (0,L),\, \tau\in (0,T) \\
u(x,0) &= f(x) & x\in (0,L)\\
u(0,\tau) = u(L,\tau) &= 0 & \tau\in (0,T).
\end{aligned}
\end{equation}
In this time-dependent differential equation, the source term on the right-hand side is zero and therefore does not require parametrization. Instead, our focus lies on the initial condition $u(x,0) = f(x)$, which characterizes the system's state at the initial time $\tau=0$. Consequently, we aim to parametrize the function $f$. 

We are interested in the well-known solution of the heat equation
\begin{align}\label{eq:heat_exact}
u(x,\tau) = \sum_{\ell=1}^\infty a_\ell \sin \left(\frac{\ell \pi x}{L}\right) \exp \left( \frac{-\ell^2\pi^2\alpha^2 \tau}{L^2} \right) & & x \in [0,L], \tau\in [0,T],
\end{align}
which can be derived exactly for the initial condition
\begin{align}\label{eq:heat_initial}
u(x,0) = f(x) = \sum_{\ell=1}^\infty a_\ell \sin \frac{\ell \pi x}{L} & & x \in [0,L]
\end{align}
with $a_\ell \in \C, \ell \in \N$, using Fourier's approach. 

While this solution holds for arbitrary $\tau \geq 0$, we set the final time $T=1$. Further, we set the length $L=1$ and the diffusivity constant $\alpha = 0.25$. Due to the time-dependence of the differential equation \eqref{eq:heat_1d}, the function space $\U$ is a little more complicated than in the previous examples. We need to ensure spatial regularity, meaning that $u$ should be square-integrable in time and satisfy $u(\tau) \in H_0^1(\Omega)$ for almost every time $\tau$, since the term $\partial_{xx}u$ and homogeneous boundary conditions are present. Additionally, we require a notion of time regularity: the weak time derivative of $u$ should exist and take values in the dual space $H^{-1}(\Omega)$. A formal and compact way to express this space uses so-called Bochner spaces (see, e.g., \cite[Sec.~5.9.2]{Evans10} or \cite{BaDu24}), which we omit here for simplicity. Accordingly, we have $\F =  L_2(\Omega)$ to ensure that the initial state $u(x,0)$ is square-integrable over the spatial domain $\Omega$. Due to the time dependence, the solution operator we are analyzing this time is of the form $\G: \F \times [0,T] \rightarrow \U$, cf. Remark \ref{rem:parameters}.

We parametrize the function $f$ by truncating the sum \eqref{eq:heat_initial} to $n$ terms. Similar to Section \ref{subsubsec:poisson_1d_fc}, we restrict the coefficients $a_\ell \in [-1,1]$ for all $\ell = 1,\ldots,n$ and transform both the spatial and time variable by $\mathcal{T} x = \frac12(x+1)$ and $\mathcal{T} \tau = \frac12(\tau+1)$. The differential equation is solved using SciPy's function \texttt{solve\textunderscore ivp}. In particular, we use the Radau IIA implicit Runge-Kutta method of order five, suitable for stiff differential equations with a desired relative tolerance of $10^{-8}$ and absolute tolerance of $10^{-10}$. We choose the number of coefficients $n=9$ to end up with an $11$-dimensional approximation problem, the sparsity $s=1000$ and the extension $N=64$ for our algorithm. Note that using \texttt{solve\textunderscore ivp} is a very costly approach for the solution of the heat equation, needing roughly $20$ seconds for a single sample of the solution. Therefore, the execution of our full algorithm took around $32$ hours.

\begin{figure}[t]
\centering
\begin{tikzpicture}
\begin{axis}[
	boxplot/draw direction=x,
	xmode=log,
	xmin=0.000009,
	xmax=0.0002,
	ytick={1},
	yticklabels={},
	xlabel={approximation error $\text{err}(\bolda)$},
	width=0.95\textwidth,
    height=0.20\textwidth,
	]
	\pgfplotstableread[col sep=comma]{5_err.csv}\datatable;
	\addplot+[
  		boxplot prepared from table={
    			table=\datatable,
    			row=0,
    			lower whisker=lw,
    			upper whisker=uw,
    			lower quartile=lq,
    			upper quartile=uq,
    			median=med
  		}, my boxplot style,
    ] coordinates {};
    \end{axis}
\end{tikzpicture}
\caption{The relative approximation error $\text{err}(\bolda)$ for $10000$ randomly drawn $\bolda$ for the heat equation example. The box-and-whisker plots show the median, the first and the second quartile as well as the maximal and minimal error observed.}
\label{fig:heat_err}
\end{figure}
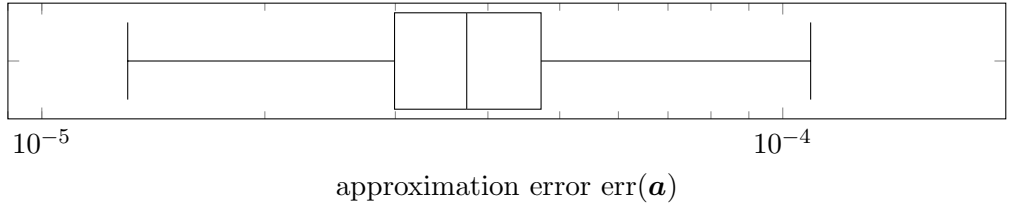

Figure \ref{fig:heat_err} shows the relative approximation error $\text{err}(\bolda)$ computed for $10000$ randomly drawn coefficients $\bolda$. The error is computed similarly to \eqref{eq:err} for $100$ equidistant nodes in space and time each and using the exact solution \eqref{eq:heat_exact} as comparison. The average error size is less than $10^{-4}$ for this example. Note, that the function \texttt{solve\textunderscore ivp} itself reaches an accuracy of around $10^{-6}$ with the given parameters, which can be seen as the noise on the data samples we are giving to our dimension-incremental method here.

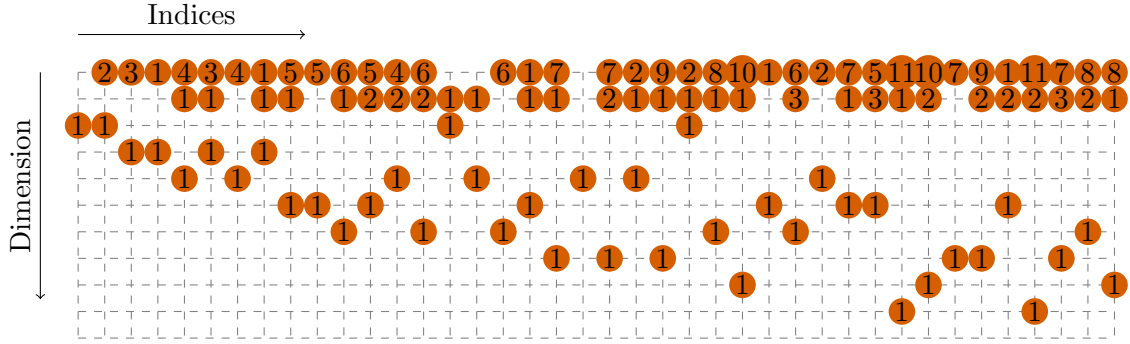
\begin{figure}[t]
\centering
\begin{tikzpicture}
\pgfplotstableread{ind.dat}\mytable
\foreach \i in {0,...,39} {
    \draw[dashed, thin, gray, opacity=0.5] (\i*10 pt,0) -- (\i*10 pt,-100 pt);
}
\foreach \j in {0,...,10} {
    \draw[dashed, thin, gray, opacity=0.5] (0,-\j*10 pt) -- (390 pt,-\j*10 pt);
}
\foreach \i in {200,...,239}{
    \foreach \j in {0,...,10}{
    \pgfplotstablegetelem{\i}{\j}\of\mytable
    \ifnum\pgfplotsretval=0\relax\else
    \node[circle, minimum size=10pt, inner sep=0pt, fill=color6, opacity=0.1*\pgfplotsretval, text opacity= 0.75] at (\i*10-2000 pt,-\j*10 pt) {\pgfplotsretval};
    \fi
    };
};
\draw[->]        (-0.5,0)  -- node [above,midway,rotate=90] {Dimension} (-0.5,-3);
\draw[->]        (0,0.5)   -- node [above,midway] {Indices} (3,0.5);
\end{tikzpicture}
\caption{An abstract visualization of the first 40 indices $\boldk$ detected for the heat equation example. The indices $\boldk$ are ordered by the absolute values of their corresponding approximated coefficients $\hat{u}_\boldk$ in descending order from left to right. The rows identify the $11$ dimensions corresponding to the spatial variable $x$, the time variable $\tau$ and the $n=9$ coefficients $\bolda$ used. Zeros are neglected to preserve clarity.}
\label{fig:heat_ind}
\end{figure}

The detected index set $I$ shows several interesting features this time. As can be seen for the first $40$ indices in Figure \ref{fig:heat_ind}, the dimensions corresponding to the coefficients $a_\ell$ contain exactly one non-zero entry for each index, which happens to be $1$. This stays true for all the detected indices up to some artifacts again. As in previous examples, this is due to the fact, that the coefficients $a_\ell$ appear only linearly and separated from each other in the sum in \eqref{eq:heat_exact}. Our algorithm once again captures this behavior even for this more complicated example. The size of the entries in the first dimension, so the dimension corresponding to the spatial variable $x$ seems to grow rapidly and much faster than in the second dimension, which corresponds to the time variable $\tau$. However, the largest entries in the first dimension is $24$ while the largest one in the second dimension is $22$. For increased sparsity $s$, the size of the entries in the first dimension stays the same while larger entries in the second dimension are added. This again confirms, that the time-dependent exponential term in \eqref{eq:heat_exact} is the main difficulty in this approximation problem.

\subsection{Burgers' equation}
\label{subsec:burger_eq}

Our final example will be the non-linear Burgers' equation in one dimension with homogeneous boundary conditions, i.e.
\begin{equation} \label{eq:burger_1d}
\begin{aligned} 
\partial_\tau u + u \partial_x u &= \nu \partial_{xx} u, & x \in (0,L),\, \tau\in (0,T) \\
u(x,0) &= f(x) & x\in (0,L)\\
u(0,\tau) = u(L,\tau) &= 0 & \tau\in (0,T).
\end{aligned}
\end{equation}
We will use the viscosity $\nu = 0.05$ here. As in the previous example, we set the final time $T=1$ and the length $L=1$, i.e. $\Omega = (0,1)$. However, we are interested in the solution $u$ only at the final time $T=1$ this time, i.e. $u(x,1)$. This again complicates the notion of our solution space $\U$ a little, which is why we omit it here. The idea however is similar to the one described for the heat equation, adjusted accordingly for the neglection of the time variable~$\tau$.

We use the approximation via the sine expansion
\begin{align}\label{eq:burger_initial_approx}
f(x) \approx \sum_{\ell=1}^n a_\ell \sin (\ell \pi x) & & x \in [0,1]
\end{align}
for the initial condition $f$ as our parametrization. This approach is similar to the one in Section \ref{subsec:heat_eq}. Analogously, we restrict $\bolda \in [-1,1]^n$, transform $\mathcal{T} x = \frac12(x+1)$ and set the number of coefficients $n=9$, such that we end up with a $10$-dimensional approximation problem. Further, we use the sparsity $s=500$ and the extension $N=32$ here. We use the function \texttt{solve\textunderscore ivp} with the Radau method again, this time with the relative tolerance $10^{-6}$ and absolute tolerance ${10^{-8}}$. The solver itself then needs about $6$ seconds for each solution. However, the more difficult structure of the index set we are detecting due to the non-linearity of our problem leaves our full method with a runtime of roughly $114$ hours.

\begin{figure}[t]
\centering
\begin{tikzpicture}
\begin{axis}[
	boxplot/draw direction=x,
	xmode=log,
	xmin=0.001,
	xmax=4,
	ytick={1},
	yticklabels={},
	xlabel={approximation error $\text{err}(\bolda)$},
	width=0.95\textwidth,
    height=0.20\textwidth,
	]
	\pgfplotstableread[col sep=comma]{6_err.csv}\datatable;
	\addplot+[
  		boxplot prepared from table={
    			table=\datatable,
    			row=0,
    			lower whisker=lw,
    			upper whisker=uw,
    			lower quartile=lq,
    			upper quartile=uq,
    			median=med
  		}, my boxplot style,
    ] coordinates {};
    \end{axis}
\end{tikzpicture}
\caption{The relative approximation error $\text{err}(\bolda)$ for $10000$ randomly drawn $\bolda$ for the Burgers' equation example. The box-and-whisker plots show the median, the first and the second quartile as well as the maximal and minimal error observed.}
\label{fig:burger_err}
\end{figure}
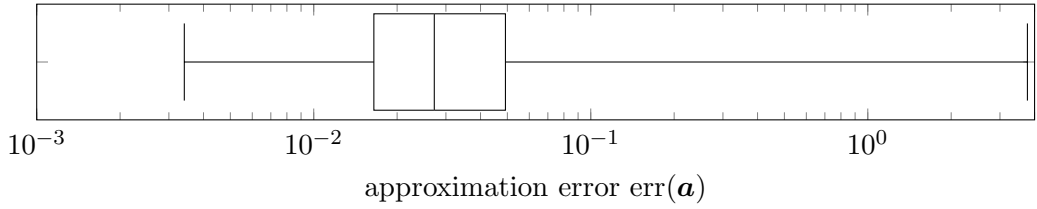

In Figure~\ref{fig:burger_err} we directly observe how the relative approximation error $\text{err}(\bolda)$ is way larger than for the heat equation in Section~\ref{subsec:heat_eq}. The non-linearity of Burgers' equation increases the difficulty of this approximation problem significantly, even when neglecting the time variable $\tau$ here. However, the error is in the order of magnitude of $10^{-2}$ most of the times, which is still reasonable given this difficulty and the smaller sparsity $s=500$ used here. Note that approximation errors $\text{err}(\bolda)$ exceeding $10^{-1}$ occur only rarely and would be classified as outliers under the common convention of using $1.5$ times the interquartile range.

Next, we are interested in an explicit solution to Burgers' equation. If we use the initial condition 
\begin{align} \label{eq:burger_initial}
f(x) = 2 \pi \nu \frac{\sin(\pi x)}{\alpha + \cos(\pi x)} & & x \in [0,L],
\end{align}
the solution $u$ is given as
\begin{align} \label{eq:burger_solution}
u(x,\tau) = 2 \pi \nu \frac{\sin(\pi x) \exp(-\pi^2 \nu \tau)}{\alpha + \cos(\pi x)\exp(-\pi^2 \nu \tau)} & & x \in [0,L], \tau \in [0,T].
\end{align}
Further details and other explicit solutions to Burgers' equation can be found in \cite{BoAwLaMuAs18}. We set $\alpha = 2$ here and approximate \eqref{eq:burger_initial} by an $n=9$-term sine expansion of the form \eqref{eq:burger_initial_approx}, which yields a relative error of less than $10^{-5}$. Note that the corresponding coefficients $a_\ell$ have absolute values of less than $0.2$, matching our inital restriction $a_\ell \in [-1,1]$.

\begin{figure}[t]
\centering
\begin{subfigure}[c]{\textwidth}
    \centering
    \begin{tikzpicture}
    \begin{axis}[
        xmin = 0,
        xmax = 1,
        legend style={at={(0.5,0.1)},anchor=south},
        width = \textwidth,
        height = 0.3\textwidth,
    ]
    \addplot[solid, color=color5, line width=1.5pt] table[x index=0, y index=2] {6_ex_err.txt};
    \addplot[dashed, color=color6, line width=1.5pt] table[x index=0, y index=1] {6_ex_err.txt};
    \legend{Exact solution, Approximate solution}
    \end{axis}
    \end{tikzpicture}
\end{subfigure}
~
\begin{subfigure}[c]{\textwidth}
    \centering
    \begin{tikzpicture}
    \begin{axis}[
        xmin = 0,
        xmax = 1,
        ylabel={Absolute error},
        width = \textwidth,
        height = 0.3\textwidth,
    ]
    \addplot[solid, color=black] table[x index=0, y index=3] {6_ex_err.txt};
    \end{axis}
    \end{tikzpicture}
\end{subfigure}
\caption{Comparison between the exact solution \eqref{eq:burger_solution} and the approximate solution at $\tau=1$ and the corresponding (absolute) pointwise approximation error.}
\label{fig:burger_explicit}
\end{figure}
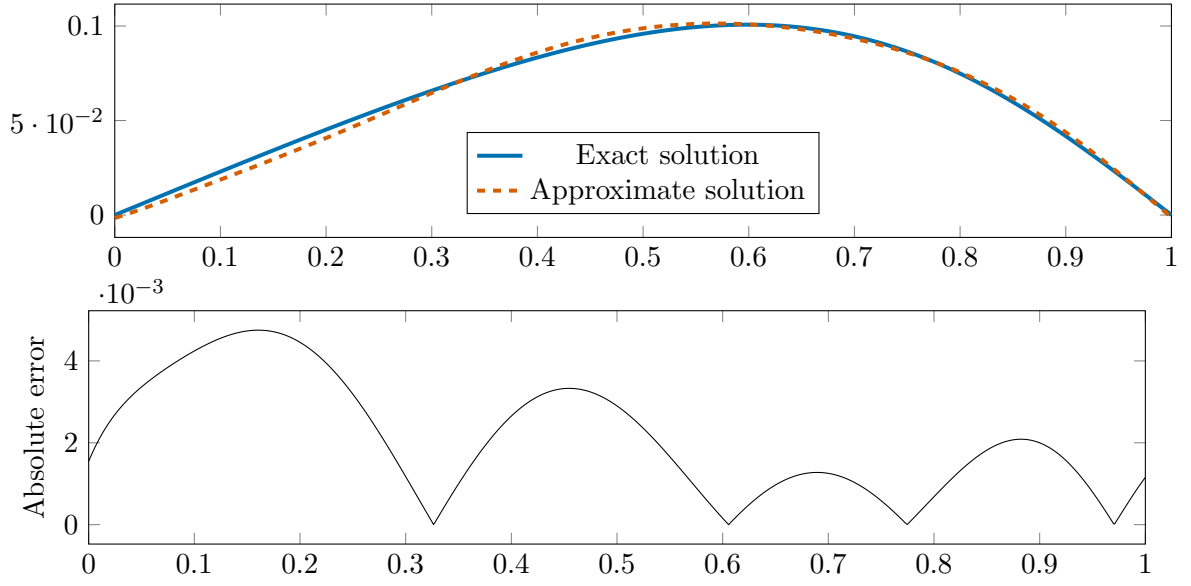

Although we are starting with just an approximation of the true initial condition \eqref{eq:burger_initial}, we end up with a reasonable solution as shown in Figure~\ref{fig:burger_explicit}. We notice the smooth oscilliations of our approximation around the true solution, which could of course be further reduced by increasing the sparsity $s$ and the extension $N$.

\begin{figure}[t]
\centering
\begin{subfigure}[c]{\textwidth}
\begin{tikzpicture}
\pgfplotstableread{ind.dat}\mytable
\foreach \i in {0,...,39} {
    \draw[dashed, thin, gray, opacity=0.5] (\i*10 pt,0) -- (\i*10 pt,-90 pt);
}
\foreach \j in {0,...,9} {
    \draw[dashed, thin, gray, opacity=0.5] (0,-\j*10 pt) -- (390 pt,-\j*10 pt);
}
\foreach \i in {240,...,279}{
    \foreach \j in {0,...,9}{
    \pgfplotstablegetelem{\i}{\j}\of\mytable
    \ifnum\pgfplotsretval=0\relax\else
    \node[circle, minimum size=10pt, inner sep=0pt, fill=color6, opacity=0.\pgfplotsretval, text opacity= 0.75] at (\i*10-2400 pt,-\j*10 pt) {\pgfplotsretval};
    \fi
    };
};
\draw[->]        (-0.5,0)  -- node [above,midway,rotate=90] {Dimension} (-0.5,-3);
\draw[->]        (0,0.5)   -- node [above,midway] {Indices} (3,0.5);
\end{tikzpicture}
\caption{The detected indices from number $1$ to number $40$.}
\end{subfigure}
~
\begin{subfigure}[c]{\textwidth}
\centering
\begin{tikzpicture}
\pgfplotstableread{ind.dat}\mytable
\foreach \i in {0,...,39} {
    \draw[dashed, thin, gray, opacity=0.5] (\i*10 pt,0) -- (\i*10 pt,-90 pt);
}
\foreach \j in {0,...,9} {
    \draw[dashed, thin, gray, opacity=0.5] (0,-\j*10 pt) -- (390 pt,-\j*10 pt);
}
\foreach \i in {280,...,319}{
    \foreach \j in {0,...,9}{
    \pgfplotstablegetelem{\i}{\j}\of\mytable
    \ifnum\pgfplotsretval=0\relax\else
    \node[circle, minimum size=10pt, inner sep=0pt, fill=color6, opacity=0.1*\pgfplotsretval, text opacity= 0.75] at (\i*10-2800 pt,-\j*10 pt) {\pgfplotsretval};
    \fi
    };
};
\draw[->]        (-0.5,0)  -- node [above,midway,rotate=90] {Dimension} (-0.5,-3);
\draw[->]        (0,0.5)   -- node [above,midway] {Indices} (3,0.5);
\end{tikzpicture}
\caption{The detected indices from number $201$ to number $240$.}
\end{subfigure}
\caption{Abstract visualizations of $40$ detected indices $\boldk$ (from left to right) for Example \ref{subsec:burger_eq}. The indices $\boldk$ are sorted in descending order according to the size of the corresponding approximated coefficient $\hat{u}_{\boldk}$. The rows identify the $10$ dimensions corresponding to the variables $x$ and $a_{1},\ldots ,a_{9}$ from top to bottom in this order. Zeros are neglected to preserve clarity.}
\label{fig:burger_ind}
\end{figure}
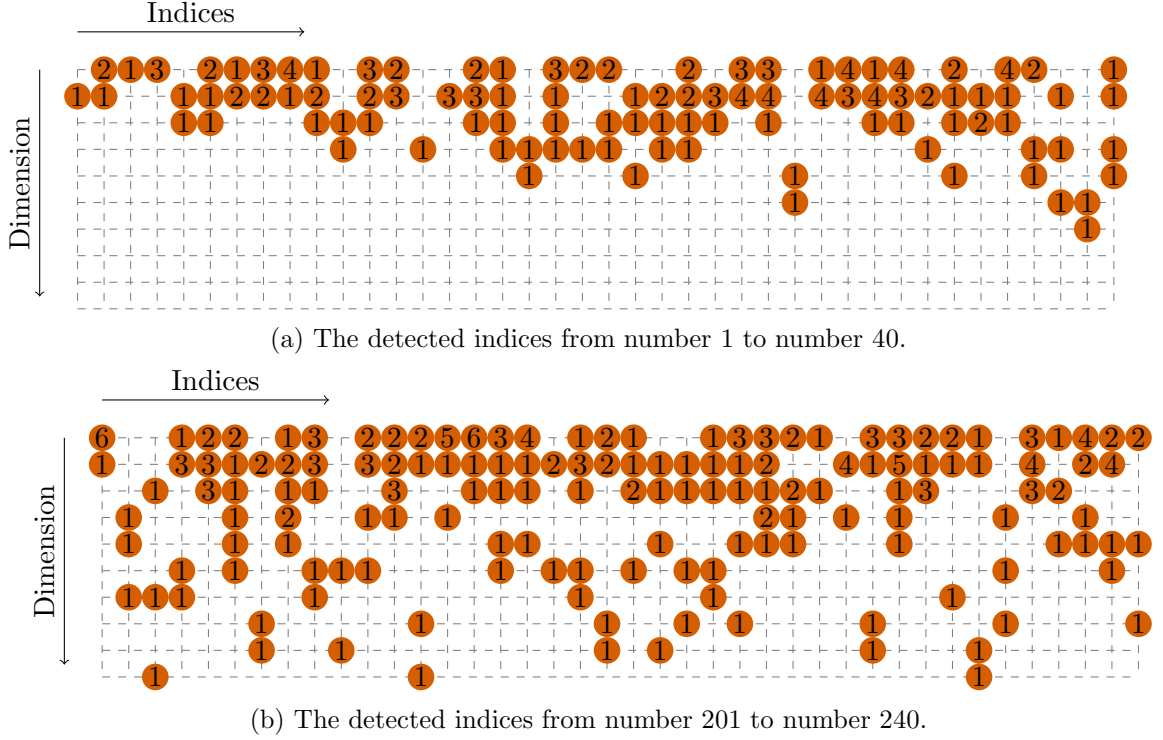

Finally, Figure~\ref{fig:burger_ind} shows some of the detected indices for this example. We observe similar behavior in the first dimension corresponding to the spatial dimension $x$ as in previous examples. However, for the remaining dimensions corresponding to the coefficients $a_\ell$, two main differences can be seen: First, the entries in these dimension are not only $0$ and $1$, i.e., the solution does not depend only linearly on the coefficients $a_\ell$. Second, there are several interactions between these dimensions (and the spatial dimension $x$) noticeable. While we saw both effects already in the example in Section~\ref{subsubsec:poisson_1d_bs}, they are much more prominent here, i.e. occuring way earlier (already in the first $10$ indices) and stronger (up to six non-zero entries).

\section{Conclusion}
\label{sec:conclusion}

We presented an adaptive approach that uses the dimension-incremental algorithm from \cite{KaPoTa23} in combination with classical differential equation solvers like the FEM in order to approximate solution operators of differential equations. We transformed the problem of operator learning for differential equations by parametrizing e.g. the source function $f$, which led us to a high-dimensional approximation problem for a function with an unknown structure. Algorithm \ref{alg:main} detects a reasonable index set $I$ by using samples of the solution $u$ computed by the differential equation solver mentioned above. This index set $I$ not only allows a good approximation of the respective solution $u$ for any source function $f$ suitable to this parametrization, but also gives us important information about the structure of the solution and its dependence on the spatial variable $\boldx$, the discretization parameters of the source function $f$ and possible other variables and parameters such as the time $t$. Such information can then be used to manually generalize the index set $I$ to even higher dimensions that arise when refining the resolution of the source function $f$. 

We have studied the behavior of our proposed methods on several examples. These numerical tests yielded reasonable approximations to the solutions of the PDEs. Especially for the easier examples, the structure of the obtained index sets $I$ matched our general expectations and (if available) the structure of the underlying analytical solution. Our brief test of generalization of the index set $I$ to even higher dimensions for the one-dimensional Poisson equation also showed promising results. The more advanced examples demonstrated the applicability of our approach to more difficult settings and problems. Our numerical tests are available at \url{https://github.com/fabiantaubert/nabopb} together with the dimension-incremental algorithm itself.

Overall, the presented algorithm performed satisfactorily and provided useful details about the structure of the solutions to the differential equations. Thus, while the field of operator learning is strongly dominated by machine learning algorithms such as PINNs, more classical approaches such as our proposed method can open new perspectives, especially to overcome still existing drawbacks of neural networks like the lack of interpretability. The 




\begin{thebibliography}{10}

\bibitem{BaDu24}
F.~Bartel and D.~Dũng.
\newblock Sampling recovery in {B}ochner spaces and applications to parametric
  {PDE}s with log-normal random inputs.
\newblock {\em arXiv:2409.05050}, 2024.

\bibitem{BaTa23}
F.~Bartel and F.~Taubert.
\newblock Nonlinear approximation with subsampled rank-1 lattices.
\newblock {\em Fourteenth International Conference on Sampling Theory and
  Applications}, 2023.

\bibitem{BlEr21}
J.~Blechschmidt and O.~G. Ernst.
\newblock Three ways to solve partial differential equations with neural
  networks — a review.
\newblock {\em GAMM-Mitteilungen}, 44(2):e202100006, 2021.

\bibitem{BoAwLaMuAs18}
M.~P. Bonkile, A.~Awasthi, C.~Lakshmi, V.~Mukundan, and V.~S. Aswin.
\newblock A systematic literature review of {B}urgers' equation with recent
  advances.
\newblock {\em Pramana - Journal of Physics}, 90(6):69, 2018.

\bibitem{BoTo24}
N.~Boullé and A.~Townsend.
\newblock Chapter 3 - a mathematical guide to operator learning.
\newblock In S.~Mishra and A.~Townsend, editors, {\em Numerical Analysis Meets
  Machine Learning}, volume~25 of {\em Handbook of Numerical Analysis}, pages
  83--125. Elsevier, 2024.

\bibitem{CaHe23}
J.~E.~S. Cardona and M.~Hecht.
\newblock Learning partial differential equations by spectral approximates of
  general {S}obolev spaces.
\newblock {\em arXiv:2301.04887}, 2023.

\bibitem{CoKuNuSu16}
R.~Cools, F.~Y. Kuo, D.~Nuyens, and G.~Suryanarayana.
\newblock Tent-transformed lattice rules for integration and approximation of
  multivariate non-periodic functions.
\newblock {\em J. Complexity}, 36:166--181, 2016.

\bibitem{DeTeGi23}
N.~Demo, M.~Tezzele, and G.~Rozza.
\newblock A {DeepONet} multi-fidelity approach for residual learning in reduced
  order modeling.
\newblock {\em Adv. Model. and Simul. in Eng. Sci.}, 10(12), 2023.

\bibitem{EiGiSchwZa14}
M.~Eigel, C.~J. Gittelson, C.~Schwab, and E.~Zander.
\newblock Adaptive stochastic {G}alerkin {FEM}.
\newblock {\em Comput. Methods Appl. Mech. Engrg.}, 270:247--269, 2014.

\bibitem{Evans10}
L.~C. Evans.
\newblock {\em Partial Differential Equations}, volume~19 of {\em Graduate
  Studies in Mathematics}.
\newblock American Mathematical Society, 2nd edition, 2010.

\bibitem{XiYuWa23}
X.~Feng, Y.~Qian, and W.~Shen.
\newblock {MC-Nonlocal-PINNs}: Handling nonlocal operators in {PINNs} via
  {M}onte {C}arlo sampling.
\newblock {\em Numer. Math. Theor. Meth. Appl.}, 16(3):769--791, 2023.

\bibitem{GoYiYuKa22}
S.~Goswami, M.~Yin, Y.~Yu, and G.~E. Karniadakis.
\newblock A physics-informed variational {DeepONet} for predicting crack path
  in quasi-brittle materials.
\newblock {\em Comput. Methods Appl. Mech. Engrg.}, 391:114587, 2022.

\bibitem{GrHeKl24}
V.~Grimm, A.~Heinlein, and A.~Klawonn.
\newblock A short note on solving partial differential equations using
  convolutional neural networks.
\newblock In {\em Domain Decomposition Methods in Science and Engineering
  XXVII}, pages 3--14, Cham, 2024. Springer Nature Switzerland.

\bibitem{GrIw24}
C.~Gross and M.~Iwen.
\newblock Sparse spectral methods for solving high-dimensional and multiscale
  elliptic {PDE}s.
\newblock {\em Found. Comput. Math.}, 2024.

\bibitem{GrKoLaSch23}
T.~G. Grossmann, U.~J. Komorowska, J.~Latz, and C.-B. Schönlieb.
\newblock Can physics-informed neural networks beat the finite element method?
\newblock {\em IMA Journal of Applied Mathematics}, 89(1):143--174, 05 2024.

\bibitem{HaSchShToTrWa22}
A.~Hashemi, H.~Schaeffer, R.~Shi, U.~Topcu, G.~Tran, and R.~Ward.
\newblock Generalization bounds for sparse random feature expansions.
\newblock {\em Applied and Computational Harmonic Analysis}, 62:310--330, 2023.

\bibitem{HeSchwaZe24}
L.~Herrmann, C.~Schwab, and J.~Zech.
\newblock Neural and spectral operator surrogates: unified construction and
  expression rate bounds.
\newblock {\em Adv. Comput. Math. (accepted)}, 2024.

\bibitem{HuNg22}
M.~Hutzenthaler and T.~A. Nguyen.
\newblock Multilevel {P}icard approximations of high-dimensional semilinear
  partial differential equations with locally monotone coefficient functions.
\newblock {\em Appl. Numer. Math.}, 181:151--175, 2022.

\bibitem{JiMeLu22}
P.~Jin, S.~Meng, and L.~Lu.
\newblock {MIONet}: Learning multiple-input operators via tensor product.
\newblock {\em SIAM J. Sci. Comput.}, 44(6):A3490--A3514, 2022.

\bibitem{Kae16}
L.~K{\"{a}}mmerer.
\newblock Multiple rank-1 lattices as sampling schemes for multivariate
  trigonometric polynomials.
\newblock {\em J. Fourier Anal. Appl.}, 24:17--44, 2018.

\bibitem{Kae25}
L.~K\"ammerer.
\newblock An efficient spatial discretization of spans of multivariate
  {C}hebyshev polynomials.
\newblock {\em Appl. Comput. Harmon. Anal.}, 77:101761, 2025.

\bibitem{KaPoTa22}
L.~K\"ammerer, D.~Potts, and F.~Taubert.
\newblock The uniform sparse {FFT} with application to {PDEs} with random
  coefficients.
\newblock {\em Sampl. Theory Signal Proces. Data Anal.}, 20(19), 2022.

\bibitem{KaPoTa23}
L.~K\"ammerer, D.~Potts, and F.~Taubert.
\newblock Nonlinear approximation in bounded orthonormal product bases.
\newblock {\em Sampl. Theory Signal Proces. Data Anal.}, 21(19), 2023.

\bibitem{KaPoVo17}
L.~K\"{a}mmerer, D.~Potts, and T.~Volkmer.
\newblock High-dimensional sparse {FFT} based on sampling along multiple rank-1
  lattices.
\newblock {\em Appl. Comput. Harmon. Anal.}, 51:225--257, 2021.

\bibitem{Ka21}
G.~E. Karniadakis, I.~G. Kevrekidis, L.~Lu, P.~Perdikaris, S.~Wang, and
  L.~Yang.
\newblock Physics-informed machine learning.
\newblock {\em Nat. Rev. Phys.}, 3:422--440, 2021.

\bibitem{Ni23}
N.~Kovachki, Z.~Li, B.~Liu, K.~Azizzadenesheli, K.~Bhattacharya, A.~Stuart, and
  A.~Anandkumar.
\newblock Neural operator: Learning maps between function spaces with
  applications to {PDEs}.
\newblock {\em J. Mach. Learn. Res.}, 24(89):1--97, 2023.

\bibitem{KoLaStu24}
N.~B. Kovachki, S.~Lanthaler, and A.~M. Stuart.
\newblock Chapter 9 - operator learning: Algorithms and analysis.
\newblock In S.~Mishra and A.~Townsend, editors, {\em Numerical Analysis Meets
  Machine Learning}, volume~25 of {\em Handbook of Numerical Analysis}, pages
  419--467. Elsevier, 2024.

\bibitem{LaStTr24}
S.~Lanthaler, A.~M. Stuart, and M.~Trautner.
\newblock Discretization error of {F}ourier neural operators.
\newblock {\em arXiv:2405.02221}, 2024.

\bibitem{Leveque07}
R.~J. LeVeque.
\newblock {\em Finite Difference Methods for Ordinary and Partial Differential
  Equations: Steady-State and Time-Dependent Problems}.
\newblock SIAM, Philadelphia, PA, 2007.

\bibitem{Li21}
Z.~Li, N.~B. Kovachki, K.~Azizzadenesheli, B.~liu, K.~Bhattacharya, A.~Stuart,
  and A.~Anandkumar.
\newblock Fourier neural operator for parametric partial differential
  equations.
\newblock In {\em International Conference on Learning Representations}, 2021.

\bibitem{LiMoZh23}
G.~Lin, C.~Moya, and Z.~Zhang.
\newblock {B-DeepONet}: An enhanced {B}ayesian {DeepONet} for solving noisy
  parametric {PDEs} using accelerated replica exchange {SGLD}.
\newblock {\em J. Comput. Phys.}, 473:111713, 2023.

\bibitem{LiMiPeKaMi23}
L.~Lingsch, M.~Michelis, S.~M. Perera, R.~K. Katzschmann, and S.~Mishra.
\newblock Vandermonde neural operators.
\newblock {\em arXiv:2305.19663}, 2023.

\bibitem{LuZh23}
T.~Luo and Q.~Zhou.
\newblock {On Residual Minimization for PDEs: Failure of PINN, Modified
  Equation, and Implicit Bias}.
\newblock {\em arXiv 2310.18201}, 2023.

\bibitem{RaPeKa19}
M.~Raissi, P.~Perdikaris, and G.~Karniadakis.
\newblock Physics-informed neural networks: A deep learning framework for
  solving forward and inverse problems involving nonlinear partial differential
  equations.
\newblock {\em J. Comput. Phys.}, 378:686--707, 2019.

\bibitem{gitNABOPB}
F.~Taubert.
\newblock {NABOPB: Nonlinear approximation in bounded orthonormal product bases
  (Python implementation)}, 2025.
\newblock Version: v1.0.0, Contributor: L.~K\"ammerer, Link:
  doi.org/10.5281/zenodo.15437933.

\bibitem{Trefethen2000}
L.~N. Trefethen.
\newblock {\em Spectral Methods in \textsc{Matlab}}.
\newblock \textrm{SIAM}, Philadelphia, PA, USA, 2000.

\bibitem{WeHuJeKr19}
E.~Weinan, M.~Hutzenthaler, A.~Jentzen, and T.~Kruse.
\newblock On multilevel {P}icard numerical approximations for high-dimensional
  nonlinear parabolic partial differential equations and high-dimensional
  nonlinear backward stochastic differential equations.
\newblock {\em J. Sci. Comput.}, 79:1534--1571, 2019.

\bibitem{WuHa21}
C.-F. Wu and M.~Hamada.
\newblock {\em Experiments: planning, analysis and parameter design
  optimization}.
\newblock Wiley series in probability and statistics. Wiley, Hoboken, NJ, third
  edition edition, 2021.

\bibitem{zech18}
J.~Zech.
\newblock {\em Sparse-Grid Approximation of High-Dimensional Parametric
  {PDE}s}.
\newblock Doctoral thesis, ETH Zurich, 2018.

\bibitem{ZiTa05}
O.~C. Zienkiewicz and R.~L. Taylor.
\newblock {\em The Finite Element Method: Its Basis and Fundamentals}.
\newblock Elsevier, Amsterdam, 6th edition edition, 2005.

\end{thebibliography}
\end{document}